\let\myfrac=\frac
\input eplain
\let\frac=\myfrac
\input amstex
\input epsf

% Pour utiliser dvips, il faut employer :

% dvips -N0 -Z0 -K0 whatever.dvi -o whatever.ps

% Here we load the functions permitting us to use amsmath without using amsppt.

\loadeufm \loadmsam \loadmsbm
\message{symbol names}\UseAMSsymbols\message{,}

\font\myfontdefault=cmr10

% cmbx10 at 16pt%
\font\mytdmchapfont=cmb10 at 14pt
\font\mytdmheadfont=cmb10 at 10pt
\font\mytdmsubheadfont=cmr10

\magnification 1200
\newif\ifintroduction
\newif\ifmaintheorems
\newif\ifinappendices
\newif\ifundefinedreferences
\newif\ifchangedreferences
\newif\ifloadreferences
\newif\ifmakebiblio
\newif\ifmaketdm

\undefinedreferencestrue
\changedreferencesfalse

% With \loadreferencesfalse, the file makes its own counters and saves them to "references.tex". With \loadreferencestrue, the file loads the counters
% saved in "references.tex".

\loadreferencestrue
\makebibliofalse
\maketdmfalse

\def\headpenalty{-400}     % 400
\def\proclaimpenalty{-200} % 200

%%%%%%%%%%%%%%%%%%%%%%%%%%%%%%%%%%%%%%%%%%%%%%%%%%%%%%%%%%%%%%%%%%%%%%%%%%%%%%%%%%%%%%%%%%%%%%%%%%%%%%%%%%%%%%%%%%%%%%%
%
% Compteurs
%
%%%%%%%%%%%%%%%%%%%%%%%%%%%%%%%%%%%%%%%%%%%%%%%%%%%%%%%%%%%%%%%%%%%%%%%%%%%%%%%%%%%%%%%%%%%%%%%%%%%%%%%%%%%%%%%%%%%%%%%

\def\alphanum#1{\ifcase #1 _\or A\or B\or C\or D\or E\or F\or G\or H\or I\or J\or K\or L\or M\or N\or O\or P\or Q\or R\or S\or T\or U\or V\or W\or X\or Y\or Z\fi}
\def\latinnum#1{\ifcase #1 _\or I\or II\or III\or IV\or V\or VI\or VII\or VIII\or IX\or X\or XI\or XII\or XIII\or XIV\or XV\or XVI\or XVII\or XVIII\or XIX\or XX\fi}
\def\gobbleeight#1#2#3#4#5#6#7#8{}

\newwrite\references
\newwrite\tdm
\newwrite\biblio

\newcount\chapno
\newcount\headno
\newcount\subheadno
\newcount\procno
\newcount\figno
\newcount\citationno

\def\setcatcodes{%
\catcode`\!=0 \catcode`\\=11}%

\ifloadreferences
    {\catcode`\@=11 \catcode`\_=11%
    \global\def\_@citation@Alexander{1}
\global\def\_@citation@Almgren{2}
\global\def\_@citation@CaffNirSprV{3}
\global\def\_@citation@Dold{4}
\global\def\_@citation@ElwTrombI{5}
\global\def\_@citation@ElwTrombII{6}
\global\def\_@citation@EspGalMira{7}
\global\def\_@citation@RosEsp{8}
\global\def\_@citation@Fethi{9}
\global\def\_@citation@Kato{10}
\global\def\_@citation@LabB{11}
\global\def\_@citation@LabI{12}
\global\def\_@citation@PacardXu{13}
\global\def\_@citation@Nirenberg{14}
\global\def\_@citation@Pogorelov{15}
\global\def\_@citation@SmiDTI{16}
\global\def\_@citation@SchneiderI{17}
\global\def\_@citation@SchneiderII{18}
\global\def\_@citation@ShengUrbasWang{19}
\global\def\_@citation@SmiAAT{20}
\global\def\_@citation@SmiSLC{21}
\global\def\_@citation@SmiNLD{22}
\global\def\_@citation@SmiPPG{23}
\global\def\_@citation@Stein{24}
\global\def\_@citation@Tromb{25}
\global\def\_@citation@White{26}
\global\def\_@citation@Ye{27}
\global\def\_@proc@ThmYeGeneralised{I}
\global\def\_@proc@ThmSignature{II}
\global\def\_@proc@ThmUniqueness{III}
\global\def\_@head@HeadPreliminaries{2}
\global\def\_@subhead@SubHeadCurvatures{2.1}
\global\def\_@head@HeadPreliminaries{2}
\global\def\_@proc@LemmaHessianOfRestrictionGeneral{2.1}
\global\def\_@proc@LemmaHessianOfRestriction{2.2}
\global\def\_@proc@LemmaHarmonicDecomposition{2.3}
\global\def\_@proc@LemmaFeynmanIntegrals{2.4}
\global\def\_@proc@ThmAlexander{2.5}
\global\def\_@proc@LemmaAlexandrovAlexander{2.6}
\global\def\_@proc@ThmAlexandrov{2.7}
\global\def\_@proc@PropLocalGeodesicProperty{2.8}
\global\def\_@proc@LemmaLocalGeodesicAlternative{2.10}
\global\def\_@head@HeadExistence{3}
\global\def\_@subhead@SubHeadExistence{3.1}
\global\def\_@proc@ThmYe{3.1}
\global\def\_@proc@LemmaExpansionOfMetric{3.2}
\global\def\_@proc@LemmaExpansionOfCurvature{3.3}
\global\def\_@proc@LemmaExpansionOfSecondFF{3.4}
\global\def\_@proc@LemmaExpansionOfShapeOperator{3.5}
\global\def\_@proc@LemmaExpansionOfGamma{3.6}
\global\def\_@proc@PropPropertiesOfK{3.7}
\global\def\_@proc@PropExpansionOfCurvature{3.8}
\global\def\_@proc@PropDiffOfMetric{3.9}
\global\def\_@proc@PropDiffOfSecondFF{3.10}
\global\def\_@proc@PropDiffOfShapeOperator{3.11}
\global\def\_@proc@PropDiffOfCurvature{3.12}
\global\def\_@proc@PropSecondDiffOfCurvature{3.13}
\global\def\_@proc@PropExpansionOfPertShapeOperator{3.14}
\global\def\_@proc@PropExpansionOfPertCurvature{3.15}
\global\def\_@proc@PropFourthOrderCoeffOfCurvature{3.16}
\global\def\_@proc@PropPertNormal{3.17}
\global\def\_@proc@PropProjectionOfCurvature{3.18}
\global\def\_@proc@PropExpansionOfFunction{3.19}
\global\def\_@proc@PropFirstAndSecondOrderTerms{3.20}
\global\def\_@proc@PropDerivativeOfFunctional{3.21}
\global\def\_@proc@PropYeGeneralised{3.22}
\global\def\_@proc@CorUniqueness{3.23}
\global\def\_@subhead@SubHeadPerturbationTheory{4.2}
\global\def\_@proc@PropExpansionOfRestriction{4.1}
\global\def\_@proc@PropZeroethAndFirstOrderCoefficients{4.2}
\global\def\_@proc@PropSecondOrderSoln{4.3}
\global\def\_@proc@PropSecondOrderOperator{4.4}
\global\def\_@proc@PropNormalVariation{4.5}
\global\def\_@proc@PropNormalVariationOfYeSpheres{4.6}
\global\def\_@proc@PropNormalComponentOfVariation{4.7}
\global\def\_@proc@PropFourthOrderOperator{4.8}
\global\def\_@proc@PropFourthOrderCoeffOfJacobiOperator{4.9}
\global\def\_@proc@ThmWeakerUniqueness{5.1}
\global\def\_@proc@PropFirstRefinement{5.3}
\global\def\_@proc@PropSecondRefinement{5.4}
\global\def\_@proc@PropThirdRefinement{5.5}
\global\def\_@proc@PropFromWeakToStrongUniqueness{5.6}
    }%
\else
    \openout\references=references.tex
\fi

\newcount\newchapflag % Flag used for checking that we have a new chapter.
\newcount\showpagenumflag % Flag used for showing page numbers.

\global\chapno = -1 % We set it like this so that the introduction does not show up a number.
\global\citationno=0
\global\headno = 0
\global\subheadno = 0
\global\procno = 0
\global\figno = 0

\def\resetcounters{%
\global\headno = 0%
\global\subheadno = 0%
\global\procno = 0%
\global\figno = 0%
}

\global\newchapflag=0 % default - false.
\global\showpagenumflag=0 % default - false.

\def\chinfo{\ifinappendices\alphanum\chapno\else\the\chapno\fi}%
\def\headinfo{\ifinappendices\alphanum\headno\else\the\headno\fi}%
\def\subheadinfo{\headinfo.\the\subheadno}%      {\chinfo.\the\headno.\the\subheadno}%
\def\procinfo{%
\ifintroduction%
    \ifmaintheorems%
        \latinnum\procno%
    \else%
        \alphanum\procno%
    \fi%
\else%
    \headinfo.\the\procno%
\fi
}%
\def\figinfo{\the\figno}%
%{\headinfo.\the\figno}
%{\chinfo.\the\headno.\the\figno}
\def\citationinfo{\the\citationno}%
\def\nextheadno{\global\advance\headno by 1 \global\subheadno = 0 \global\procno = 0}% \headinfo\ - }
\def\nextsubheadno{\global\advance\subheadno by 1}% \subheadinfo\ - }
\def\nextprocno{\global\advance\procno by 1 \procinfo}
\def\nextfigno{\global\advance\figno by 1 \figinfo}

{\global\let\noe=\noexpand%
%
% Ici, je definit les macros qui me permettent de contruire des compteurs. Je vais changer provisoirement les catcodes afin de pouvoir bien ecrire le fichier
% contenant les references.
%
\catcode`\@=11%
\catcode`\_=11%
\setcatcodes%
!global!def!_@@internal@@makeref#1{%
!global!expandafter!def!csname #1ref!endcsname##1{%
!csname _@#1@##1!endcsname%
!expandafter!ifx!csname _@#1@##1!endcsname!relax%
    !write16{#1 ##1 not defined - run saving references}%
    !undefinedreferencestrue%
!fi}}%
!global!def!_@@internal@@makelabel#1{%
!global!expandafter!def!csname #1label!endcsname##1{%
!edef!temptoken{!csname #1info!endcsname}%
!ifloadreferences%
    !expandafter!ifx!csname _@#1@##1!endcsname!relax%
        !write16{#1 ##1 not hitherto defined - rerun saving references}%
        !changedreferencestrue%
    !else%
        !expandafter!ifx!csname _@#1@##1!endcsname!temptoken%
        !else
            !write16{#1 ##1 reference has changed - rerun saving references}%
            !changedreferencestrue%
        !fi%
    !fi%
!else%
    !expandafter!edef!csname _@#1@##1!endcsname{!temptoken}%
    !edef!textoutput{!write!references{\global\def\_@#1@##1{!temptoken}}}%
    !textoutput%
!fi}}%
!global!def!makecounter#1{!_@@internal@@makelabel{#1}!_@@internal@@makeref{#1}}%
!unsetcatcodes%
}
\makecounter{ch}%
\makecounter{head}%
\makecounter{subhead}%
\makecounter{proc}%
\makecounter{fig}%
\makecounter{citation}%
\def\newref#1#2{%
\def\temptext{#2}%
\edef\bibliotextoutput{\expandafter\gobbleeight\meaning\temptext}%
\global\advance\citationno by 1\citationlabel{#1}%
\ifmakebiblio%
    \edef\fileoutput{\write\biblio{\noindent\hbox to 0pt{\hss$[\the\citationno]$}\hskip 0.2em\bibliotextoutput\medskip}}%
    \fileoutput%
\fi}%
\def\cite#1{%
$[\citationref{#1}]$%
\ifmakebiblio%
    \edef\fileoutput{\write\biblio{#1}}%
    \fileoutput%
\fi%
}%
%
%%%%%%%%%%%%%%%%%%%%%%%%%%%%%%%%%%%%%%%%%%%%%%%%%%%%%%%%%%%%%%%%%%%%%%%%%%%%%%%%%%%%%%%%%%%%%%%%%%%%%%%%%%%%%%%%%%%%%%%
%
% Mise en page
%
%%%%%%%%%%%%%%%%%%%%%%%%%%%%%%%%%%%%%%%%%%%%%%%%%%%%%%%%%%%%%%%%%%%%%%%%%%%%%%%%%%%%%%%%%%%%%%%%%%%%%%%%%%%%%%%%%%%%%%%

\let\mypar=\par

% Here we define right justification.

\def\raggedleft{\leftskip=0pt plus 1fil \parfillskip=0pt}

% Lettrine - grace a Serroul

\font\lettrinefont=cmr10 at 28pt
\def\lettrine #1[#2][#3]#4%
{\hangafter -#1 \hangindent #2
\noindent\hskip -#2 \vtop to 0pt{
\kern #3 \hbox to #2 {\lettrinefont #4\hss}\vss}}

\font\mylettrinefont=cmr10 at 28pt
\def\mylettrine #1[#2][#3][#4]#5%
{\hangafter -#1 \hangindent #2
\noindent\hskip -#2 \vtop to 0pt{
\kern #3 \hbox to #2 {\mylettrinefont #5\hss}\vss}}

% Here we define the default page headers and footers.

\edef\Pagetitle={Blank}

\headline={\hfil\Pagetitle\hfil}
%{
%\ifnum\showpagenumflag=0
%    \hfil
%\else
%    \ifnum\newchapflag=0
%        \ifodd\pageno
%            \hfil\tenrm{\sl\Pagetitle}\ -\ \the\pageno
%        \else
%            \tenrm\hss\kern-1.2mm\the\pageno\ -\ {\sl\Pagetitle}\hfill
%        \fi
%    \else
%        \hfil
%    \fi
%\fi
%}

\footline={\hfil\myfontdefault\folio\hfil}
%{
%\ifnum\showpagenumflag=0
%    \hfil
%\else
%    \ifnum\newchapflag=0
%        \hfil
%    \else
%        \global\newchapflag=0
%        \ifodd\pageno
%            \hfil\tenrm\the\pageno
%        \else
%            \tenrm\the\pageno\hfil
%        \fi
%    \fi
%\fi
%}

\def\nextoddpage
{
\newpage%
\ifodd\pageno%
\else%
    \global\showpagenumflag = 0%
    \null%
    \vfil%
    \eject%
    \global\showpagenumflag = 1%
\fi%
}

% Ici on definit les titres des chapitres

\def\newchap#1#2%
{%
%
% First we reset the counters.
%
\global\advance\chapno by 1%
\resetcounters%
%
% Next we move on to the next page.
%
\newpage%
\ifodd\pageno%
\else%
    \global\showpagenumflag = 0%
    \null%
    \vfil%
    \eject%
    \global\showpagenumflag = 1%
\fi%
\global\newchapflag = 1%
\global\showpagenumflag = 1%
%
% Now we write the chapter heading.
%
{\font\chapfontA=cmsl10 at 30pt%
\font\chapfontB=cmsl10 at 25pt%
\null\vskip 5cm%
{\chapfontA\raggedleft\hfil%
{%
\ifnum\chapno=0
    \phantom{%
    \ifinappendices%
        Annexe \alphanum\chapno%
    \else%
        \the\chapno%
    \fi}%
\else%
    \ifinappendices%
        Annexe \alphanum\chapno%
    \else%
        \the\chapno%
    \fi%
\fi%
}%
\par}%
\vskip 2cm%
{\chapfontB\raggedleft%
\lineskiplimit=0pt%
\lineskip=0.8ex%
\hfil #1\par}%
\vskip 2cm%
}%
\edef\Pagetitle{#2}%
%
% Finally, we write the information into the tdm file.
%
\ifmaketdm%
    \def\temp{#2}%
    \def\tempbis{\nobreak}%
    \edef\chaptitle{\expandafter\gobbleeight\meaning\temp}%
    \edef\mynobreak{\expandafter\gobbleeight\meaning\tempbis}%
    \edef\textoutput{\write\tdm{\bigskip{\noexpand\mytdmchapfont\noindent\chinfo\ - \chaptitle\hfill\noexpand\folio}\par\mynobreak}}%
\fi%
\textoutput%
}

% Ici on definit les sous titres.

\def\newhead#1%
{%
\ifhmode%
    \mypar%
\fi%
\ifnum\headno=0%
\ifinappendices
    \nobreak\vskip -\lastskip%
    \nobreak\vskip .5cm%
\fi
\else%
    \nobreak\vskip -\lastskip%
    \nobreak\vskip .5cm%
\fi%
\nextheadno%
\ifmaketdm%
    \def\temp{#1}%
    \edef\sectiontitle{\expandafter\gobbleeight\meaning\temp}%
    \edef\textoutput{\write\tdm{\noindent{\noexpand\mytdmheadfont\quad\headinfo\ - \sectiontitle\hfill\noexpand\folio}\par}}%
    \textoutput%
\fi%
\font\headfontA=cmbx10 at 14pt%
{\headfontA\noindent\headinfo\ - #1.\hfil}%
\nobreak\vskip .5cm%
}%

% Ici on definit les sous-sous titres.

\def\newsubhead#1%
{%
\ifhmode%
    \mypar%
\fi%
\ifnum\subheadno=0%
\else%
    \penalty\headpenalty\vskip .4cm%
\fi%
\nextsubheadno%
\ifmaketdm%
    \def\temp{#1}%
    \edef\subsectiontitle{\expandafter\gobbleeight\meaning\temp}%
    \edef\textoutput{\write\tdm{\noindent{\noexpand\mytdmsubheadfont\quad\quad\subheadinfo\ - \subsectiontitle\hfill\noexpand\folio}\par}}%
    \textoutput%
\fi%
\font\subheadfontA=cmsl10 at 12pt% cmbxsl10 at 12pt% cmbxti110 at 12pt%
{\subheadfontA\noindent\subheadinfo\ #1.\hfil}%
\nobreak\vskip .25cm %
}%

%%%%%%%%%%%%%%%%%%%%%%%%%%%%%%%%%%%%%%%%%%%%%%%%%%%%%%%%%%%%%%%%%%%%%%%%%%%%%%%%%%%%%%%%%%%%%%%%%%%%%%%%%%%%%%%%%%%%%%%
%
% Commandes et Symboles
%
%%%%%%%%%%%%%%%%%%%%%%%%%%%%%%%%%%%%%%%%%%%%%%%%%%%%%%%%%%%%%%%%%%%%%%%%%%%%%%%%%%%%%%%%%%%%%%%%%%%%%%%%%%%%%%%%%%%%%%%

% Ici on definit la famille mathroman...

\font\mathromanten=cmr10
\font\mathromanseven=cmr7
\font\mathromanfive=cmr5
\newfam\mathromanfam
\textfont\mathromanfam=\mathromanten
\scriptfont\mathromanfam=\mathromanseven
\scriptscriptfont\mathromanfam=\mathromanfive
\def\mathroman{\fam\mathromanfam}

% Ici on definit la famille mathsf qui n'est pas definie dans Plain.

\font\sf=cmss12

\font\sansseriften=cmss10
\font\sansserifseven=cmss7
\font\sansseriffive=cmss5
\newfam\sansseriffam
\textfont\sansseriffam=\sansseriften
\scriptfont\sansseriffam=\sansserifseven
\scriptscriptfont\sansseriffam=\sansseriffive
\def\mathsf{\fam\sansseriffam}

% Ici on definit la famille mathbf qui n'est pas definie dans Plain.

\font\bftwelve=cmb12

\font\boldten=cmb10
\font\boldseven=cmb7
\font\boldfive=cmb5
\newfam\mathboldfam
\textfont\mathboldfam=\boldten
\scriptfont\mathboldfam=\boldseven
\scriptscriptfont\mathboldfam=\boldfive
\def\mathbf{\fam\mathboldfam}

% Ici on definit \mathi et \mathj. En fait, cela n'est pas necessaire, car ces symboles existent deja sous les noms respectivement de \imath et \jmath. Mais,
% ca reste une exercice interessant.

\font\mycmmiten=cmmi10
\font\mycmmiseven=cmmi7
\font\mycmmifive=cmmi5
\newfam\mycmmifam
\textfont\mycmmifam=\mycmmiten
\scriptfont\mycmmifam=\mycmmiseven
\scriptscriptfont\mycmmifam=\mycmmifive

\def\hexa#1{\ifcase #1 0\or 1\or 2\or 3\or 4\or 5\or 6\or 7\or 8\or 9\or A\or B\or C\or D\or E\or F\fi}
\mathchardef\mathi="7\hexa\mycmmifam7B
\mathchardef\mathj="7\hexa\mycmmifam7C

% Ici je definit \mybeth, \mygimmel et \mydaleth

\font\mymsbmten=msbm10 at 8pt%8.5pt
\font\mymsbmseven=msbm7 at 5.6pt%6pt
\font\mymsbmfive=msbm5 at 4pt%4.25pt
\newfam\mymsbmfam
\textfont\mymsbmfam=\mymsbmten
\scriptfont\mymsbmfam=\mymsbmseven
\scriptscriptfont\mymsbmfam=\mymsbmfive

\mathchardef\mybeth="7\hexa\mymsbmfam69
\mathchardef\mygimmel="7\hexa\mymsbmfam6A
\mathchardef\mydaleth="7\hexa\mymsbmfam6B

% Ici je definit la fonction qui insere des figures.

\def\placelabel[#1][#2]#3{{%
\setbox10=\hbox{\raise #2cm \hbox{\hskip #1cm #3}}%
\ht10=0pt%
\dp10=0pt%
\wd10=0pt%
\box10}}%

% Ici je definit les fonctions \proclaim et \endproclaim. Les premieres versions sont plus jolies mais je ne sais pas comment les couper verticalement, et il y
% a alors des problemes de mise en page. Les deuxiemes versions sont plus simples, mais il n'y a pas de problem de mise en page.

\newif\ifinproclaim%
\global\inproclaimfalse%
\def\proclaim#1{%
\medskip%
%
% J'insere ici un bgroup pour que les variables que je definis ici soit uniquement locales.
%
\bgroup%
\inproclaimtrue%
\setbox10=\vbox\bgroup\leftskip=0.8em\noindent{\bftwelve #1}\sf%
}

\def\endproclaim{%
\egroup%
\setbox11=\vtop{\noindent\vrule height \ht10 depth \dp10 width 0.1em}%
\wd11=0pt%
\setbox12=\hbox{\copy11\kern 0.3em\copy11\kern 0.3em}%
\wd12=0pt%
\setbox13=\hbox{\noindent\box12\box10}%
\noindent\unhbox13%
\egroup%
\medskip\ignorespaces%
}

\def\proclaim#1{%
\medskip%
\bgroup%
\inproclaimtrue%
\noindent{\bftwelve #1}%
\nobreak\medskip%
\sf%
}

\def\endproclaim{%
\mypar\egroup\penalty\proclaimpenalty\medskip\ignorespaces%
}

\def\noskipproclaim#1{%
\medskip%
\bgroup%
\inproclaimtrue%
\noindent{\bf #1}\nobreak\sl%
}

\def\endnoskipproclaim{%
\mypar\egroup\penalty\proclaimpenalty\medskip\ignorespaces%
}

% Ici on definit quelques autres commandes dont on aura besoin, en particulier, le colon francais.

\def\ninn{{n\in\Bbb{N}}}

\def\proof{{\noindent\bf Proof:\ }}

\def\remark{{\noindent\sl Remark:\ }}

\def\mlim{\mathop{{\mathroman Lim}}}

\def\msf#1{{\mathsf #1}}

\def\qed{~$\square$}

\def\minter{\mathop{\cap}}
\def\myitem#1{%
%\ifinproclaim%
%    \item{#1}%
%\else%
    \noindent\hbox to .5cm{\hfill#1\hss}%\kern 0.1em}%
}
%\fi}%

\catcode`\@=11
\def\Eqalign#1{\null\,\vcenter{\openup\jot\m@th\ialign{%
\strut\hfil$\displaystyle{##}$&$\displaystyle{{}##}$\hfil%
&&\quad\strut\hfil$\displaystyle{##}$&$\displaystyle{{}##}$%
\hfil\crcr #1\crcr}}\,}
\catcode`\@=12

\def\makeop#1{%
\global\expandafter\def\csname op#1\endcsname{{\mathroman #1}}}%

\def\makeopsmall#1{%
\global\expandafter\def\csname op#1\endcsname{{\mathroman{\lowercase{#1}}}}}%

\makeopsmall{ArcTan}%
\makeopsmall{ArcCos}%
\makeop{Arg}%
\makeop{Det}%
\makeop{Log}%
\makeop{Re}%
\makeop{Im}%
\makeop{Dim}%
\makeopsmall{Tan}%
\makeop{Ker}%
\makeopsmall{Cos}%
\makeopsmall{Sin}%
\makeop{Exp}%
\makeopsmall{Tanh}%
\makeop{Tr}%
\makeop{End}%
\makeop{Long}%
\makeop{Ch}%
\makeop{Exp}%
\makeop{Eval}%
\makeop{Lift}%
\makeop{Int}%
\makeop{Ext}%
\makeop{Aire}%
\makeop{Im}%
\makeop{mult}%
\makeop{Conf}%
\makeop{Crit}
\makeop{Exp}%
\makeopsmall{Mod}%
\makeop{Inj}%
\makeop{Log}%
\makeop{Sgn}%
\makeop{Ext}%
\makeop{Int}%
\makeop{Ln}%
\makeop{Dist}%
\makeop{Aut}%
\makeop{Id}%
\makeop{GL}%
\makeop{SO}%
\makeop{Homeo}%
\makeop{Imm}%
\makeop{Vol}%
\makeop{Ric}%
\makeop{Hess}%
\makeop{Euc}%
\makeop{Scal}%
\makeop{Isom}%
\makeop{Max}%
\makeop{SW}%
\makeop{SL}%
\makeop{Long}%
\makeop{Alex}%
\makeop{Fix}%
\makeop{Wind}%
\makeop{Diag}%
\makeop{dVol}%
\makeop{Coker}%
\makeop{Symm}%
\makeop{Ad}%
\makeop{Diam}%
\makeop{loc}%
\makeopsmall{Sinh}%
\makeopsmall{Cosh}%
\makeop{Len}%
\makeop{Mult}%
\makeop{Length}%
\makeop{Conv}%
\makeop{Min}%
\makeop{Area}%
\font\mycirclefont=cmsy7
\def\textcircle{{\raise 0.3ex \hbox{\mycirclefont\char'015}}}

\let\emph=\bf

\hyphenation{quasi-con-formal}

%%%%%%%%%%%%%%%%%%%%%%%%%%%%%%%%%%%%%%%%%%%%%%%%%%%%%%%%%%%%%%%%%%%%%%%%%%%%%%%%%%%%%%%%%%%%%%%%%%%%%%%%%%%%%%%%%%%%%%%
%
% Citations
%
%%%%%%%%%%%%%%%%%%%%%%%%%%%%%%%%%%%%%%%%%%%%%%%%%%%%%%%%%%%%%%%%%%%%%%%%%%%%%%%%%%%%%%%%%%%%%%%%%%%%%%%%%%%%%%%%%%%%%%%

\ifmakebiblio%
    \openout\biblio=biblio.tex %
    {%
        \edef\fileoutput{\write\biblio{\bgroup\leftskip=2em}}%
        \fileoutput
    }%
\fi%

\newref{Alexander}{Alexander S., Locally convex hypersurfaces of negatively curved spaces, {\sl Proc. Amer. Math. Soc.} {\bf 64} (1977), no. 2, 321--325}
\newref{Almgren}{Almgren F. J. Jr., {\sl Plateau's problem: An invitation to varifold geometry},.W. A. Benjamin, Inc., New York-Amsterdam, (1966)}
%\newref{Berger}{Berger M., {\sl A panoramic view of Riemannian geometry}, Springer-Verlag, Berlin, (2003)}
%\newref{BoehmWilking}{B\"ohm C., Wilking B., Manifolds with positive curvature operators are space forms, {\sl Ann. of Math. (2)} {\bf 167} (2008), no. 3, 1079--1097}
\newref{CaffNirSprV}{Caffarelli L., Nirenberg L., Spruck J., Nonlinear second-order elliptic equations. V. The Dirichlet problem for Weingarten hypersurfaces, {\sl Comm. Pure Appl. Math.} {\bf 41} (1988), no. 1, 47--70}
\newref{Dold}{Dold A., {\sl Lectures on Algebraic Topology}, Classics in Mathematics, Springer-Verlag, Berlin, Heidelberg, (1995)}
\newref{ElwTrombI}{Elworthy K. D., Tromba A. J., Differential structures and Fredholm maps on Banach manifolds, in {\sl Global Analysis (Proc. Sympos. Pure Math.)}, Vol. {\bf XV}, Berkeley, Calif., (1968), 45--94 Amer. Math. Soc., Providence, R.I.}
\newref{ElwTrombII}{Elworthy K. D., Tromba A. J., Degree theory on Banach manifolds, in {\sl Nonlinear Functional Analysis (Proc. Sympos. Pure Math.)}, Vol. {\bf XVIII}, Part 1, Chicago, Ill., (1968), 86--94, Amer. Math. Soc., Providence, R.I.}
\newref{EspGalMira}{Espinar J. M., G\'alvez J. A., Mira P., Hypersurfaces in $\Bbb H^{n+1}$ and conformally invariant equations: the generalized Christoffel and Nirenberg problems, {\sl J. Eur. Math. Soc.}, {\bf 11}, no. 4, (2009), 903--939}
\newref{RosEsp}{Espinar J., Rosenberg H., When strictly locally convex hypersurfaces are embedded, arXiv:1003.0101}
\newref{Fethi}{Fethi M., Constant $k$-curvature hypersurfaces in Riemannian manifolds, {\sl Differential Geom. Appl.} {\bf 28} (2010), no. 1, 1--11}
\newref{Kato}{Kato T., {\sl Perturbation theory for linear operators}, Die Grundlehren der mathematischen Wissenschaften, {\bf 132}, Springer-Verlag, New York, (1966)}
\newref{LabB}{Labourie F., Probl\`eme de Minkowski et surfaces \`a courbure constante dans les vari\'et\'es hyperboliques, {\sl Bull. Soc. Math. France}, {\bf 119}, no. 3,(1991), 307--325}
\newref{LabI}{Labourie F., Immersions isométriques elliptiques et courbes pseudoholomorphes, {\sl J. Diff. Geom.} \bf{30}, (1989), 395-44}
\newref{PacardXu}{Pacard F., Xu X., Constant mean curvature spheres in Riemannian manifolds, {\sl Manuscripta Math.} {\bf 128} (2009), no. 3, 275--295}
%\newref{MicallefMoore}{Micallef M. J., Moore J. D., Minimal two-spheres and the topology of manifolds with positive curvature on totally isotropic two-planes, {\sl Ann. of Math.} (2), \bf{127}, (1988), no. 1, 199--227}
\newref{Nirenberg}{Nirenberg L., The Weyl and Minkowski problem in differential geometry in the large, {\sl Comm. Pure Appl. Math.}, {\bf 6}, no. 3, (1953), pp. 337--394}
\newref{Pogorelov}{Pogorelov A.V., Extrinsic geometry of convex surfaces, {\sl Israel program for scientific translation}, Jerusalem, (1973)}
\newref{SmiDTI}{Rosenberg H., Smith G., Degree Theory of Immersed Hypersurfaces, arXiv:1010.1879}
\newref{SchneiderI}{Schneider M., Alexandrov Embedded closed magnetic geodesics on $S^2$, arXiv:0903.1128}
\newref{SchneiderII}{Schneider M., Closed magnetic geodesics on closed hyperbolic Riemann Surfaces, arXiv:1009.1723}
%\newref{Schwarz}{Schwarz M., {\sl Morse Homology}, Progress in Mathematics, {\bf 111}, Birkh\"auser Verlag, Basel, %Boston, Berlin, (1993)}
\newref{ShengUrbasWang}{Sheng W., Urbas J., Wang X., Interior curvature bounds for a class of curvature equations. (English summary), {\sl Duke Math. J.} {\bf 123} (2004), no. 2, 235--264}
\newref{SmiAAT}{Smith G., An Arzela-Ascoli Theorem for Immersed Submanifolds, {\sl Ann. Fac. Sci. Toulouse Math.} {\bf 16} (2007), no. 4, 817--866}
\newref{SmiSLC}{Smith G., Special Lagrangian curvature, to appear in {\sl Math. Annalen}}
%\newref{SmiCGC}{Smith G., Constant Gaussian Curvature hypersurfaces in Hadamard manifolds, arXiv:0908.3590}
\newref{SmiNLD}{Smith G., The non-linear Dirichlet problem in Hadamard manifolds, arXiv:0908.3590}
\newref{SmiPPG}{Smith G., The Plateau problem for general curvature functions, arXiv:1008.3545}
\newref{Stein}{Stein E., {\sl Singular integrals and differentiability properties of functions}, Princeton University Press, (1970)}
\newref{Tromb}{Tromba A. J., The Euler characteristic of vector fields on Banach manifolds and a globalization of Leray-Schauder degree, {\sl Adv. in Math.}, {\bf 28}, (1978), no. 2, 148--173}
\newref{White}{White B., The space of $m$-dimensional surfaces that are stationary for a parametric elliptic functional, {\sl Indiana Univ. Math. J.}, {\bf 36}, (no. 3), (1987), 567--602}
\newref{Ye}{Ye R., Foliation by constant mean curvature spheres, {\sl Pacific J. Math.} {bf 147} (1991), no. 2, 381--396}

\ifmakebiblio%
    {\edef\fileoutput{\write\biblio{\egroup}}%
    \fileoutput}%
\fi%

%%%%%%%%%%%%%%%%%%%%%%%%%%%%%%%%%%%%%%%%%%%%%%%%%%%%%%%%%%%%%%%%%%%%%%%%%%%%%%%%%%%%%%%%%%%%%%%%%%%%%%%%%%%%%%%%%%%%%%%
%
% The Paper
%
%%%%%%%%%%%%%%%%%%%%%%%%%%%%%%%%%%%%%%%%%%%%%%%%%%%%%%%%%%%%%%%%%%%%%%%%%%%%%%%%%%%%%%%%%%%%%%%%%%%%%%%%%%%%%%%%%%%%%%%
%
\document
\myfontdefault
\global\chapno=1
\global\showpagenumflag=1
\def\Pagetitle{}
\null
\vfill
\def\centre{\rightskip=0pt plus 1fil \leftskip=0pt plus 1fil \spaceskip=.3333em \xspaceskip=.5em \parfillskip=0em \parindent=0em}%
\def\textmonth#1{\ifcase#1\or January\or Febuary\or March\or April\or May\or June\or July\or August\or September\or October\or November\or December\fi}
\font\abstracttitlefont=cmr10 at 14pt
{\abstracttitlefont\centre Constant curvature hyperspheres and the Euler Characteristic\par}
\bigskip
{\centre Graham Smith\par}
\bigskip
%{\centre 9 December 2009\par}
{\centre \the\day\ \textmonth\month\ \the\year\par}
\bigskip
{\centre Centre de Recerca Matem\`atica,\par
Campus de Bellaterra, Edifici C,\par
08193 Bellaterra,\par
Barcelona,\par
SPAIN\par}
\bigskip
\noindent{\emph Abstract:\ }We show how in many cases the algebraic number of immersed hyperspheres of constant (and prescribed) curvature may be related to the Euler Characteristic of the ambient space.
%
% 58E12 : Global analysis, analysis on manifolds; applications to minimal surfaces.
% 35J60 : Nonlinear elliptic equations
% 53A10 : Minimal surfaces, surfaces with prescribed mean curvature
% 53C21 : Methods of Riemannian geometry, including PDE methods; curvature restrictions
% 53C42 : Immersions (minimal, prescribed curvature, tight, etc.)
%
\par
\vfill
\nextoddpage
\global\pageno=1
\def\Pagetitle{\sl Constant curvature hyperspheres and the Euler Characteristic}
\newhead{Introduction}
\introductiontrue
\newsubhead{Background}
\noindent Minkowski type problems, which are now classical in the study of Riemannian submanifolds, ask for the construction of closed, immersed hypersurfaces of constant curvature subject to topological or geometric restrictions. There are many known techniques for constructing solutions, including Geometric Measure Theory, polyhedral approximation, the continuity method, curvature flow methods, and even the Perron Method (c.f. \cite{Almgren}, \cite{EspGalMira}, \cite{LabB}, \cite{Nirenberg}, \cite{Pogorelov} for a brief and very incomplete list of examples).
\medskip
\noindent Closely related to the question of existence is the more refined question of exactly how many solutions there actually are. Thus, in our forthcoming work, \cite{SmiDTI}, written in collaboration with Harold Rosenberg, we conjecture for example, that in the case where the ambient manifold is the $3$-dimensional sphere, under certain conditions on the Riemannian metric, there exist at least two immersed (even embedded) spheres of constant mean or constant extrinsic curvature. A simpler but equally interesting problem is then to count, not the number of solutions, but rather the difference between the numbers of solutions of two different types, and to this end we construct in \cite{SmiDTI} a degree theory for closed, immersed hypersurfaces of constant curvature, which extends to general notions of curvature the degree theory for constant mean curvature hypersurfaces developed by Brian White in \cite{White}, and which is closely related to the PDE degree theories of Elworthy and Tromba (c.f. \cite{ElwTrombI}, \cite{ElwTrombII} and \cite{Tromb}), as well as the degree theory of magnetic geodesics developed by Schneider in \cite{SchneiderI}.
\medskip
\noindent The real interest of such a degree theory lies in its relationships to other known topological or geometric invariants. A few results are already known. Thus, it is relatively easy to show that when the ambient manifold is diffeomorphic to a standard sphere but with non-trivial metric, under suitable curvature assumptions the algebraic number of immersed hyperspheres of prescribed curvature is generically equal to $-1$ times the Euler Characteristic of the sphere: in other words, $-2$ when the ambient sphere is even dimensional, and $0$ otherwise. Similarly, in \cite{SchneiderI} and \cite{SchneiderII}, Schneider shows that when the ambient space, $\Sigma$, is $2$-dimensional, with curvature bounded below by $-K^2$, then for generic $\kappa\in C^\infty(\Sigma)$ such that $\kappa>K$, the algebraic number of Alexandrov embedded closed curves of prescribed geodesic curvature equal to $\kappa$ is also equal to $-1$ times the Euler Characteristic of $\Sigma$.
\medskip
\noindent One is thus naturally led to suspect that this result generalises in some sense to all compact manifolds in all dimensions, but to exactly what extent, and in exactly what way presents an intriguing problem, and it as a partial response which forms the content of the current paper, as well as our forthcoming work, \cite{SmiDTI}: for a compact Riemannian manifold, $M$, we embed $\Bbb{R}$ in $C^\infty(M)$ in the obvious manner, and then, in a wide variety of cases, there exists $K>0$, which is defined by the geometry of $M$, and a neighbourhood, $\Omega$ of $]K,+\infty[$ in $C^\infty(M)$ such that for generic $\kappa\in\Omega$, the algebraic number of immersed hyperspheres in $M$ of prescribed curvature equal to $\kappa$ is itself equal to $-1$ times the Euler Characteristic of that manifold.
\medskip
\noindent Proving this result is divided into three stages: the first, to develop an appropriate degree theory in a general setting for immersed hypersurfaces of prescribed curvature; the second, to prove for any given notion of curvature the compactness properties which allow this degree theory to be applied (recall that topological degree is only defined for proper maps. Our current setting is no different in this respect, and some notion of compactness is therefore required); and the third, to show that this degree has the desired value. It is the final part which we treat in the current work. The first two parts will be studied in \cite{SmiDTI}, which will then allow us to prove results concerning the algebraic number of locally strictly convex spheres of prescribed mean, extrinsic and special Lagrangian curvature in a large class of ambient manifolds.
\newsubhead{Main Results}
\maintheoremstrue
\noindent We thus show that, in many circumstances, the algebraic number of locally strictly convex immersions of prescribed curvature (which is well defined by \cite{SmiDTI}) is equal to $-1$ times the Euler Characteristic of the ambient manifold. However, the precise statement of this result depends on the compactness properties of the curvature notion in question, for which there is still no general theory. We thus state our main results in terms of the following three theorems which, in rough terms, show that modulo the supplementary hypotheses depending on the curvature notion in question, when $k>0$ is large:
\medskip
\myitem{(i)} each critical point of the scalar curvature function generates an immersed hypersphere of constant curvature equal to $k$;
\medskip
\myitem{(ii)} the signature of this immersed hypersphere (as defined in \cite{SmiDTI}) is equal to $(-1)^n$ times the signature of the corresponding critical point of the scalar curvature function, where $(n+1)$ is the dimension of the ambient manifold; and
\medskip
\myitem{(iii)} under suitable geometric assumptions, these are the only immersed hyperspheres of constant curvature equal to $k$.
\medskip
\noindent Thus, for the notions of curvature that interest us classical index theory tells us that, grosso-modo, for sufficiently large values of $k$, the algebraic number of immersed hypersurfaces of constant curvature equal to $k$ is equal to $(-1)$ times the Euler Characteristic of the ambient manifold, and since by Theorem $1.2$ of \cite{SmiDTI} the degree is constant, the result follows. There is a clear analogy here with the behaviour of better-known theories of elliptic submanifolds (such as pseudo-holomorphic curves). Indeed, when $K$ is special Lagrangian curvature (resp. $2$-dimensional extrinsic curvature), hypersurfaces of constant $K$-curvature are in bijective correspondence with certain families of special Legendrian submanifolds of (resp. pseudo-holomorphic curves in) the unitary bundle of the ambient manifold (which we recall naturally carries a contact structure). Letting the curvature tend to infinity in the base is equivalent to taking the adiabatic limit in the total space.
\medskip
\noindent The first result is already proven in the mean curvature case by Ye in \cite{Ye}, and readily extends to general notions of curvature, with one important limitation for our purposes, being that, wheras the family of Morse functions is a generic subset of the space of smooth functions over any manifold, it is less obvious that the family of metrics whose scalar curvature function is a Morse function is generic. We readily circumvent this by considering small perturbatons of the scalar curvature function and looking for hypersurfaces not of constant, but of prescribed curvature. We obtain:
\proclaim{Theorem \nextprocno, {\bf Existence - Short Version}}
\noindent Let $R$ be the scalar curvature function of $M$, and let $f\in C^\infty(M)$ be a smooth function such that $R_f:=R + 2(n+3)/(n+1)f$ is a Morse Function. If $p\in M$ is a critical point of $R_f$ then there exists a smooth family, $(p_t,\varphi_t)_{t\in[0,\epsilon[}$, of centered spheres converging to $p$ such that, for all $t\in[0,\epsilon[$, the $K$-curvature of $(p_t,\varphi_t)$ is prescribed by the function $t^{-1}(1+t^2 f)$.
\endproclaim
\proclabel{ThmYeGeneralised}
\remark The concept of centered spheres converging to a point is explained in Section \subheadref{SubHeadExistence}.
\medskip
\noindent In the sequel, the spheres constructed by Theorem \procref{ThmYeGeneralised} about the critical points of $R_f$ will be referred to as {\bf Ye's Spheres}. The most original and most technical part of this paper lies in proving the second result, where we determine their signatures. We obtain:
\proclaim{Theorem \nextprocno, {\bf Signature}}
\noindent Let $M:=M^{n+1}$ be an $(n+1)$-dimensional Riemannian manifold. Let $R$ be the scalar curvature function of $M$ and let $f\in C^\infty(M)$ be a smooth function such that $R_f:=R + 2(n+3)/(n+1)f$ is a Morse Function. Let $p\in M$ be a critical point of $R_f$ and let $(p_t,\varphi_t)_{t\in[0,\epsilon[}$ be the family of Ye's spheres about $p$. Then, for sufficiently small $t$, $\Sigma_t:=(t,p_t,\varphi_t)$ is non-degenerate, and its signature is given by:
$$
\opSgn(\Sigma_t) = -(-1)^{n+1}\opSgn(\opHess(R_f)(p)).
$$
\endproclaim
\proclabel{ThmSignature}
\noindent Finally, to show that Ye's spheres are the only spheres that appear when the curvature is sufficently large, it is necessary to obtain bounds for their diameters, for which there is currently no general result. In the current context we therefore limit ourselves to the following case:
\proclaim{Theorem \nextprocno, {\bf Uniqueness}}
\noindent Suppose $K=H$ is mean curvature, or $K$ is convex and satisfies the Unbounded Growth Axiom (Axiom $(vii)$), then, for all $f\in C^\infty(M)$ and for all $C>0$, there exists $B>0$ such that if $H\geqslant B$, and if $\Sigma$ is a locally strictly convex immersed sphere in $M$ such that:
\medskip
\myitem{(i)} $\Sigma$ is of prescribed $K$-curvature equal to $H(1+H^{-2}f)$; and
\medskip
\myitem{(ii)} $\opDiam_M(\Sigma)\leqslant CH^{-1}$.
\medskip
\noindent Then $\Sigma$ is one Ye's spheres.
\endproclaim
\proclabel{ThmUniqueness}
\remark A statement of the Unbounded Growth Axiom is given in Section \subheadref{SubHeadCurvatures}. Importantly, this axiom is satisfied in particular, by extrinsic curvature.
\newsubhead{Summary and Acknowledgements}
\maintheoremsfalse
\global\procno=0
\noindent This paper is structured as follows:
\medskip
\myitem{(i)} in Section $2$, we review the necessary background material, introducing curvature functions, drawing our inspiration from the work, \cite{CaffNirSprV}, of Caffarelli, Nirenberg and Spruck; the basic theory of spherical harmonics; basic properties of locally strictly convex immersions; and the Geodesic Boundary Property of locally strictly convex immersions, which is important for the application of the regularity results of \cite{ShengUrbasWang};
\medskip
\myitem{(ii)} in Section $3$, we prove Theorem \procref{ThmYeGeneralised}. Although our reasoning presents little more than a mild simplification of \cite{Ye}, we provide a fairly detailled account, since the formulae obtained will be of use in the sequel;
\medskip
\myitem{(iii)} in Section $4$, we determine the signature of the Jacobi operators of the spheres constructed in Theorem \procref{ThmYeGeneralised}, thus proving Theorem \procref{ThmSignature}. This constitutes the hardest and most technical part of the paper. The signature is determined by calculating the asymptotic development of the Jacobi operator up to and including order $4$. All contributions to the signature are shown to be zero up to and including order $3$, and a geometric argument is then used to determine the fourth order contribution; and
\medskip
\myitem{(iv)} in Section $5$, we prove Theorem \procref{ThmUniqueness}, which is carried out in two stages: first we prove the general result, Theorem \procref{ThmWeakerUniqueness}, which yields uniqueness as soon as the rescaled spheres are known to converge smoothly to an embedded sphere, and Theorem \procref{ThmUniqueness} then follows by establishing various conditions under which the rescaled spheres can be shown to converge smoothly.
\medskip
\noindent This paper was written while the author was benefitting from a Marie Curie Postdoctoral Fellowship hosted by the Centre de Recerca Matem\`atica, Barcelona, Spain. The author is grateful to Harold Rosenberg for introducing him to this problem.
\newhead{Preliminaries}
\introductionfalse
\newsubhead{Curvature Functions}
\noindent We will be concerned throughout the sequel with hypersurfaces of constant and prescribed curvature for general curvature functions. We recall here the definition of curvature functions, as outlined in \cite{CaffNirSprV} (but see also \cite{SmiPPG}). Let $\Gamma_0\subseteq\Bbb{R}^n$ be the cone of vectors all of whose components are non-negative:
\headlabel{HeadPreliminaries}
\subheadlabel{SubHeadCurvatures}
$$
\Gamma = \left\{(x_1,...,x_n)\ |\ x_1,...,x_n\geqslant 0\right\}.
$$
\noindent Let $\Gamma\subseteq\Bbb{R}^n$ be a convex cone with vertex the origin. We say that $\Gamma$ is a {\bf supporting cone} for a curvature function if and only if:
\medskip
\myitem{(i)} for every permutation, $\sigma$, and for all $(x_1,...,x_n)\in\Gamma$:
$$
(x_{\sigma(1)},...,x_{\sigma(n)}) \in\Gamma;\text{\ and}
$$
\myitem{(ii)} for all $x\in\Gamma$:
$$
x + \Gamma_0 \subseteq \Gamma.
$$
\noindent A {\bf curvature function} is a pair $(K,\Gamma)$, where $\Gamma\subseteq\Bbb{R}^n$ is a supporting cone, and $K:\Gamma\rightarrow[0,\infty[$ is a non-negative valued function which is smooth over the interior of $\Gamma$ such that:
\medskip
{\noindent\bf Axiom $\mathbf{(i)}$:\ }for every permutation, $\sigma$, and for all $x_1,..,x_n\in\Gamma$:
$$
K(x_{\sigma(1)},...,x_{\sigma(n)}) = K(x_1,...,x_n);
$$
{\noindent\bf Axiom $\mathbf{(ii)}$:\ }$K$ is homogeneous of order $1$;
\medskip
{\noindent\bf Axiom $\mathbf{(iii)}$:\ }$K(1,1,...,1)=1$;
\medskip
{\noindent\bf Axiom $\mathbf{(iv)}$:\ }$K$ is strictly positive over the interior of $\Gamma$ and vanishes over $\partial\Gamma$.
\medskip
{\noindent\bf Axiom $\mathbf{(v)}$:\ }$K$ is strictly elliptic. In other words, for all $x\in\Gamma\subseteq\Bbb{R}^n$ and for all $1\leqslant i\leqslant n$:
$$
(\partial_i K)(x) > 0;\text{\ and}
$$
{\noindent\bf Axiom $\mathbf{(vi)}$:\ }$K$ is a concave function over $\Gamma\subseteq\Bbb{R}^n$.
\medskip
\noindent We say that a curvature function is {\bf convex} if and only if $\Gamma=\Gamma_0$. In the sequel, we often suppress $\Gamma$, and denote the curvature function merely by $K$. Scalar notions of curvature of hypersurfaces are generated by curvature functions in the following manner: let $M:=M^{n+1}$ be an $(n+1)$-dimensional Riemannian manifold, let $\Sigma=(S,i)$ be an immersed hypersurface in $M$, let $A$ be the shape operator of $\Sigma$, and let $K:=(K,\Gamma)$ be a curvature function. We say that $\Sigma$ is {\bf $K$-convex} if and only if the vector of eigenvalues of $A$ is an element of $\Gamma$. Since $\Gamma$ is invariant under permutation of the coordinates, this notion is well defined. We likewise define strict $K$-convexity in the obvious manner, and in the sequel, all hypersurfaces will be assumed to be strictly $K$-convex, and we shall no longer stress this property. We then define $K_\Sigma$, the {\bf $K$-curvature} of $\Sigma$, by:
$$
K_\Sigma = K(\lambda_1,...,\lambda_n),
$$
\noindent where $\lambda_1,...,\lambda_n$ are the eigenvalues of $A$. Axioms $(i)$ to $(iv)$ then become natural from a geometric perspective, whilst Axioms $(v)$ and $(vi)$ do so from an analytic perspective (see \cite{SmiPPG} for more details). The reader may check that all standard notions of curvature, such as mean curvature, extrinsic curvature, $\sigma_k$-curvature (where $\sigma_k$ is the $k$'th order symmetric polynomial), and so on, satisfy Axioms $(i)$ to $(vi)$, and that extrinsic curvature, the curvature quotients $\sigma_n/\sigma_k$ and special Lagrangian curvature (c.f. \cite{SmiSLC}) are all convex.
\medskip
\noindent In addition, for technical reasons, in this paper we will be interested in curvature functions which satisfy the following axiom, which we refer to in the sequel as the {\bf Unbounded Growth Axiom}:
\medskip
{\noindent\bf Axiom $\mathbf{(vii)}$:\ }for all $x:=(x_1,...,x_n)\in\Gamma$:
$$
\mlim_{t\rightarrow+\infty}K(x_1+t,x_2,...,x_n) = +\infty.
$$
\noindent This axiom only intervenes in the proof of Theorem \procref{ThmUniqueness}, and we leave the interested reader to verify that the following weaker version is actually sufficient for our purposes:
\medskip
{\noindent\bf Axiom $\mathbf{(vii)'}$:\ }for all $x:=(x_1,...,x_n)\in\Gamma$:
$$
\mlim_{t\rightarrow+\infty}K(x_1+t,x_2+t,x_3,...,x_n) = +\infty.
$$
\newsubhead{Spherical Harmonics}
\noindent Let $\Sigma^n$ be the $n$-dimensional unit sphere viewed as a subset of $(n+1)$-dimensional Euclidean space. Let $\opHess$ denote the Hessian operator of its Levi-Civita covariant derivative. Thus, if $X$ and $Y$ are vector fields over $\Sigma^n$ and $\alpha$ is any tensor field, then:
\headlabel{HeadPreliminaries}
$$
\opHess(\alpha)(X,Y) = \nabla_Y\nabla_X\alpha - \nabla_{\nabla_Y X}\alpha,
$$
\noindent where $\nabla$ is the Levi-Civita covariant derivative of $\Sigma^n$. Let $\opHess^0$ denote the canonical Hessian operator over $\Bbb{R}^{n+1}$. We recall the following relation for immersed hypersurfaces in general manifolds (c.f Lemma $2.8$ of \cite{SmiNLD}):
\proclaim{Lemma \nextprocno}
\noindent Let $M:=(M^{n+1},g)$ be an $(n+1)$-dimensional Riemannian manifold. Let $N:=N^n$ be an $n$-dimensional Riemannian manifold and let $i:N\rightarrow M$ be a smooth immersion. Let $\opHess^M$ and $\opHess^N$ be the Hessian operators of $(M,g)$ and $(N,i^*g)$ respectively. If $f\in C^\infty(M)$, and if $X,Y\in TN$, then:
$$
\opHess^N(f)(X,Y) = \opHess^M(f)(X,Y) - \langle\nabla f,\msf{N}\rangle II(X,Y),
$$
\noindent where $\msf{N}$ is the unit normal vector field over $N$ compatible with the orientation and $II$ is its second fundamental form.
\endproclaim
\proclabel{LemmaHessianOfRestrictionGeneral}
\noindent In particular, when $N=\Sigma^n$ is the unit sphere in Euclidean space, the second fundamental form coincides with the metric, and we immediately obtain:
\proclaim{Lemma \nextprocno}
\noindent Let $f:\Bbb{R}^{n+1}\rightarrow\Bbb{R}$ be a smooth function. If $X$, $Y$ are tangent vectors to $\Sigma^n$, then:
$$
\opHess(f)(X,Y) = \opHess^0(f)(X,Y) - \langle D f,\msf{N}\rangle\langle X,Y\rangle,
$$
\noindent where $\langle\cdot,\cdot\rangle$ is the canonical metric of $\Bbb{R}^{n+1}$, $D$ is its Levi-Civita covariant derivative, and $\msf{N}$ is the outward pointing unit normal vector field over $\Sigma^n$.
\endproclaim
\proclabel{LemmaHessianOfRestriction}
\noindent Let $\Delta$ be the Laplacian of $\Sigma_n$ acting on $C^\infty(\Sigma^n)$:
$$
\Delta f = \sum_{i=1}^n\opHess(f)(e_i,e_i),
$$
\noindent where $e_1,...,e_n$ is an orthonormal basis of $T\Sigma^n$. Let $\Bbb{N}_0$ denote the family of non-negative integers, and for all $n\in\Bbb{N}_0$, let $\Cal{H}_k\subseteq C^\infty(\Bbb{R}^{n+1})$ denote the family of harmonic, homogeneous polynomials of order $k$. We refer to the elements of $\Cal{H}_k$ as the {\bf solid spherical harmonics of order $k$}. For all $k$, let $H_k$ denote the family of restrictions of elements of $\Cal{H}_k$ to $\Sigma^n$. We refer to $H_k$ as the {\bf spherical harmonics of order $k$}. It follows from Lemma \procref{LemmaHessianOfRestriction} that, for all $k$, the elements of $H_k$ are eigenvectors of $\Delta$ with eigenvalue $\lambda_k=-k(n+k-1)$. Moreover, any eigenvector of $\Delta$ belongs to $H_k$ for some $k$, and $L^2(\Sigma^n)=\oplus_{k\in\Bbb{N}_0}H_k$ constitutes a decomposition of $L^2(\Sigma^n)$ into an orthogonal family of eigenspaces for $\Delta$ (c.f. \cite{Stein}, Ch. $III$). More precisely, if $\Cal{P}_k$ denotes the space of polynomials of order $k$ over $\Bbb{R}^{n+1}$ and if $P_k$ denotes the space of restrictions of elements of $\Cal{P}$ to $\Sigma^n$, then (c.f. \cite{Stein}, Ch. $III$):
\proclaim{Lemma \nextprocno}
\noindent For all $k\in\Bbb{N}_0$:
$$
P_{2k} = \oplus_{i=0}^k H_{2i},\qquad P_{2k+1} = \oplus_{i=0}^k H_{2i+1}.
$$
\endproclaim
\proclabel{LemmaHarmonicDecomposition}
\noindent For all $k$, let $\Pi_k:L^2(S^n)\rightarrow H_k$ be the orthogonal projection. The following result is used to obtain explicit formulae for $\Pi_1$ in the cases of interest to us:
\proclaim{Lemma \nextprocno}
\noindent Let $x_1,...,x_{n+1}$ be the coordinate functions. Let $\Omega_n$ be the volume of the $n$-dimensional sphere. For all $(i,j)$:
$$
\int_{S^n}x^ix^j\opdVol = \frac{\Omega_n}{(n+1)}\delta_{ij}.
$$
\noindent Likewise, for all $(i,j,k,l)$:
$$
\int_{S^n}x^ix^jx^kx^l\opdVol = \frac{\Omega_n}{(n+1)(n+3)}(\delta_{ij}\delta_{kl} + \delta_{ik}\delta_{jl} + \delta_{il}\delta_{jk}).
$$
\endproclaim
\proclabel{LemmaFeynmanIntegrals}
\proof Consider the $2$-tensor $A_{ij}$ given by:
$$
A_{ij} = \int_{S^n}x_ix_j\opdVol.
$$
\noindent Since $A_{ij}$ is preserved by $O(n)$, it is equal to $K\delta_{ij}$, for some constant $K$. However:
$$
\Omega_n = \int_{S^n}r^2\opdVol = (n+1)\int_{S^n}x_i^2\opdVol = (n+1)K.
$$
\noindent The first result follows. To prove the second result, we define $A_{ijkl}$ in an analogous manner. Since it is preserved by $O(n)$, it is equal to $K(\delta_{ij}\delta_{kl} + \delta_{ik}\delta_{jl} + \delta_{il}\delta_{jk})$, for some $K$, and we conclude as before.\qed
\medskip
\noindent In addition, we underline the following trivial properties which are particularly relevant to our situation:
\medskip
\myitem{(i)} $-(\Delta+n)$ is self-adjoint;
\medskip
\myitem{(ii)} all its eigenvalues are real;
\medskip
\myitem{(iii)} it has only one strictly negative real eigenvalue, $l_0:=-n$, which is of multiplicity $1$; and
\medskip
\myitem{(iv)} the eigenspace of $l_1:=0$ is $H_1$, which constitutes the restrictions to $\Sigma^n$ of the coordinate functions of $\Bbb{R}^{n+1}$, and, in particular, is of dimension $(n+1)$, and so $l_1$ has multiplicity $(n+1)$.
\medskip
\noindent We recall from \cite{SmiDTI} that if $L$ is a generalised Laplacian defined over a compact manifold, then the {\bf signature} of $L$ is defined to be equal to $(-1)^{\opmult(L)}$, where $\opmult(L)$ is equal to the number of its eigenvalues (counted with multiplicity) which are both real and strictly negative. In the sequel, we study perturbations of the operator $-\frac{1}{n}(\Delta+n)$, which is the Jacobi operator of any $K$-curvature defined over $\Sigma^n$. In order to determine the signatures of these perturbations, it suffices to study perturbations of the degenerate eigenvalue $l_1=0$.
\newsubhead{Locally Convex Immersions}
\noindent We review the properties of locally convex immersions in Riemannian manifolds. We first require the following generalised version of the Hadamard-Stoker Theorem, which we have been unable to find in the literature (but see also \cite{Alexander}):
\proclaim{Theorem \nextprocno}
\noindent Suppose that $M$ has dimension at least $3$ and that its sectional curvature is bounded above by $B$. If $i:\Sigma\rightarrow M$ is a compact, locally convex immersed hypersurface such that:
$$
\opDiam(\Sigma)<\opMin(\pi/(4\sqrt{B}),\opInj(M)/2).
$$
\noindent Then $\Sigma$ is embedded, and bounds a convex set.
\endproclaim
\proclabel{ThmAlexander}
\proof Let $\msf{N}$ be the outward pointing unit normal vector field over $\Sigma$, and define the mapping $I:\Sigma\times[0,\infty[\rightarrow M$ by:
$$
I(p,t) = \opExp(t\msf{N}(p)).
$$
\noindent Denote $R:=\opMin(\pi/(4\sqrt{B}),\opInj(M)/2)$. Using Jacobi Fields, we readily show that the restriction of $I$ to $\Sigma\times[0,2R[$ is an immersion. Choose $\delta>0$ and $p\in M$ such that $\Sigma\subseteq B:=B_{R-\delta}(p)$. Using Jacobi Fields again, we show that for all $r<2R$, and for all $q\in M$, the radial lines in $B_r(q)$ are distance minimising. This implies that $d(q,\partial B_r(q)) = r$, and so, for any $q\in B$, every geodesic, $\gamma$ leaving $q$ intersects $\partial B$ at $\gamma(t)$ for some $t\leqslant 2(R-\delta)$. We define $G:\Sigma\rightarrow\partial B$ to be the Gauss Map which sends the point $x\in\Sigma$ to the first point of intersection with $\partial B$ of the geodesic leaving $x$ in the direction of $\msf{N}(x)$.
\medskip
\noindent Using Jacobi Fields again, we show that, for all $r<2R$ and for all $q\in M$, $B_r(q)$ is strictly convex, and thus so is $B$. In particular, any geodesic inside $B$ crosses $\partial B$ transversally at the first point of contact. Since $I$ restricted to $\Sigma\times[0,2R[$ is an immersion, this transversality implies that $G$ is bilipschitz and, being a map between $2$ compact spaces, is therefore a covering map. Since $\partial B\cong S^n$ is simply connected, $G$ is a diffeomorphism, and $\Sigma$ is thus diffeomorphic to the $n$-dimensional sphere.
\medskip
\noindent For all $p\in\Sigma$, let $\Gamma_p$ be the geodesic arc leaving $p$ in the direction of $\msf{N}(p)$ and terminating in $\partial B$. Let $\Omega\subseteq\Sigma$ be the set of all points $p\in\Sigma$ such that $\Gamma_p$ only intersects $\Sigma$ at its endpoint, $p$. $\Omega$ is trivially open. Let $p\in\Sigma$ and $q\in\partial B$ be points minimising the distance between $\Sigma$ and $\partial B$. Since $\Sigma$ and $\partial B$ are locally strictly convex, the second variation formula for geodesic variations implies that $\msf{N}(p)$ points towards $q$ and so $G(p) = q$. In particular, $\Gamma_p$ only intersects $\Sigma$ at its end point, $p$, and so $p\in\Omega$. Now suppose $p\in\partial\Omega$ and let $(p_n)_\ninn\in\Omega$ be a sequence converging to $p$. $\Gamma_0 := \Gamma_{p}$ is tangential to $\Sigma$ at some point, $q$, say. For all sufficiently large $n$, there therefore exists $q_n$ close to $q$ such that some point of $\Gamma_n:=\Gamma_{p_n}$ lies above $q_n$. Since $\Sigma$ is compact and since $I$ is an immersion, elementary degree theory yields a point $p_0'\neq p_0$ such that $G(p_0')=G(p_0)$. This is absurd, and it follows that $\partial\Omega=\emptyset$. By connectedness, we deduce that $\Omega=\Sigma$, and so for all $p\in\Sigma$, $\Gamma_p$ only intersects $\Sigma$ at its endpoint, $p$. Using Degree Theory again, we deduce that $I$ is injective. In particular, $\Sigma$ is embedded, and it follows by the Jordan Sphere Theorem (c.f. \cite{Dold}) that $\Sigma$ bounds an open set which is necessarily convex. This completes the proof.\qed
\medskip
\noindent When $M$ is $2$-dimensional, the above proof fails because of non-simple connectedness of the circle which permits multiple coverings. However when the immersion is Alexandrov embedded, we recover a version of the Hadamard-Stoker Theorem:
\proclaim{Theorem \nextprocno}
\noindent Let $M$ be a $2$-dimensional Hadamard manifold, and let $i:S^1\rightarrow M$ be a locally convex, closed, immersed curve. If, in addition, $i$ is Alexandrov Embedded, then $i$ is embedded and bounds a convex set.
\endproclaim
\proclabel{LemmaAlexandrovAlexander}
\proof Taking $B=0$ in the proof of Theorem \procref{ThmAlexander}, we consider the Gauss Mapping, $G$, of $i$ from $S^1$ onto the boundary of a ball of large radius in $M$. As before, $G$ is a covering map. Moreover, since $i$ is Alexandrov embedded, we readily show that $G$ is homotopic to a map of degree $1$, and is therefore itself of degree $1$, and the rest of the proof proceeds as before.\qed
\medskip
\noindent We require the following version of the classical Alexandrov Theorem:
\proclaim{Theorem \nextprocno, {\bf Alexandrov}}
\noindent Let $K$ be any curvature function, and let $\Sigma\subseteq\Bbb{R}^{n+1}$ be a bounded, compact, locally convex immersed hypersurface of constant $K$ curvature equal to $1$. Then $\Sigma$ is a sphere of radius $1$.
\endproclaim
\proclabel{ThmAlexandrov}
\proof By Theorem \procref{ThmAlexander} with $B=0$, $\Sigma$ is embedded and bounds a convex set. The result now follows by the Alexandrov Reflection Principle.\qed
\newsubhead{The Local Geodesic Property of Convex Sets}
\noindent Let $K\subseteq\Bbb{R}^{n+1}$ be a closed, convex set we say that $K$ satisfies the {\bf Local Geodesic Property} if and only if for every point $p\in K$, there exists a geodesic segment, $\Gamma$, containing $p$ in its interior such that $\Gamma\subseteq K$. Of course, this is non-trivial only when $p\in\partial K$ is a boundary point of $K$. The following characterisation of the Local Geodesic Property is important for the application of regularity theory:
\proclaim{Proposition \nextprocno}
\noindent $K$ satisfies the Local Geodesic Property at $p\in\partial K$ if and only if for any supporting hyperplane, $H$ of $K$ at $p$, and for all sufficiently small $r>0$, the intersection $K\minter H\minter \partial B_r(p)$ is not contained in the interior of a hemisphere.
\endproclaim
\proclabel{PropLocalGeodesicProperty}
\proof If $K$ satisfies the Local Geodesic Property at $p$, and if $H$ is a supporting hyperplane of $K$ at $p$, then $H\minter K$ contains a geodesic arc passing through $p$. Thus, for all sufficiently small $r$, $K\minter H\minter\partial B_r(p)$ contains at least two antipodal points and is therefore not contained in the interior of a hemisphere. This proves the first implication.
\medskip
\noindent We prove the converse by induction. Suppose that, for every supporting hyperplane, $H$, of $K$ at $p$, and for every small $r>0$, the intersection $K\minter H\minter \partial B_r(p)$ is not contained in the interior of a hemisphere. Let $H$ be a supporting hyperplane. Trivially $K'=K\minter H\subset\Bbb{R}^n$ is convex. If $p\in K'$ is an interior point, then $K'$, and therefore $K$, trivially satisfies the Local Geodesic Property at $p$. We therefore suppose that $p\in \partial K'$ is a boundary point. Let $H'\subseteq H$ be a supporting hyperplane to $K'$ at $p$. We claim that for every small $r>0$, the intersection $K'\minter H'\minter\partial B_r(p)$ is not contained in the interior of a hemisphere. Indeed, $K\minter H\minter \partial B_r(p)$ is contained in the closed hemisphere bounded by $H'\minter\partial B_r(p)$. Suppose now that $K'\minter H'\minter \partial B_r(p)$ is contained in the interior of a hemisphere bounded by the hyperplane, $H''$, say. Then, tilting $H'$ slightly about $H''$ we obtain a hemisphere whose interior contains $K\minter H\minter\partial B_r(p)$, which is absurd. Proceeding by induction, we conclude that either $p$ satisfies the Local Geodesic Property, or there exists a convex subset $K''\subseteq\Bbb{R}^2$ which neither satisfies the Local Geodesic Property at $p$ and whose intersection with $\partial B_r(p)$ is not contained in the interior of a hemisphere for all sufficiently small $r$. This is absurd, and the result follows.\qed
\medskip
\noindent In \cite{ShengUrbasWang}, Sheng, Urbas and Wang, show that if the $C^{0,1}$ limit of a sequence of functions satisfying a certain type of non-linear elliptic PDE is a strictly convex function near that point, then it is smooth over a well defined neighbourhood of that point. The relevance of the Local Geodesic Property to this regularity theory follows from the following corollary:
\proclaim{Corollary \nextprocno}
\noindent If $K$ has non-trivial interior and does not satisfy the Local Geodesic Property at $p$, then there exists:
\medskip
\myitem{(i)} an affine hyperplane, $P$;
\medskip
\myitem{(ii)} a convex, open subset $U\subseteq P$;
\medskip
\myitem{(iii)} an open subset $V\subseteq \partial K$; and
\medskip
\myitem{(iv)} a $C^{0,1}$ function $f:U\rightarrow [0,\infty[$,
\medskip
\noindent such that:
\medskip
\myitem{(i)} $p\in V$;
\medskip
\myitem{(ii)} $V$ is the graph of $f$ over $U$; and
\medskip
\myitem{(iii)} $p$ lies strictly above $U$ (i.e. $f(p)>0$);
\medskip
\myitem{(iv)} $f$ vanishes along $\partial U$.
\endproclaim
\proof Let $H$ be a supporting hyperplane to $K$ at $p$. Since $K$ has non-trivial interior, we may suppose that $K$ is a graph over $H$ in a neighbourhood of $p$. There exists $r>0$ as small as we wish such that $K\minter H\minter \partial B_r(p)$ is contained in the interior of a hemisphere. Tilting $H$ slightly therefore yields a hyperplane, $P$, such that $K\minter P\minter \partial B_r(p)=\emptyset$. Choosing $P$ sufficiently close to $H$ and translating it slightly away from $p$, we obtain a hyperplane with the desired properties, and this completes the proof.\qed
\medskip
\noindent On the other hand, the geometric implications of the Geodesic Boundary Property are described in the following lemma:
\proclaim{Lemma \nextprocno}
\noindent Let $K\subseteq\Bbb{R}^{n+1}$ be a convex set, and let $X\subseteq\partial\Omega$ be the set of all points $p\in\partial K$ satisfying the Local Geodesic Property. Suppose that $X$ is closed. Then either:
\medskip
\myitem{(i)} $X=\emptyset$; or
\medskip
\myitem{(ii)} $X=\partial\Omega$ and $\Omega$ is unbounded.
\endproclaim
\proclabel{LemmaLocalGeodesicAlternative}
\proof Suppose that $X$ is bounded. We claim that $X=\emptyset$. Indeed, choose $p\in\Bbb{R}^{n+1}$ and let $d_p$ be the distance in $\Bbb{R}^{n+1}$ to $p$. Then, since $X$ is closed, and therefore compact, there exists $q\in X$ maximising $d_p$. However, by definition, there exists a geodesic segment, $\Gamma\in\Bbb{R}^{n+1}$ such that $q$ lies in the interior of $\Gamma$ and $\Gamma\subseteq X$. It follows that $d_p$ is maximised along $\Gamma$ at $q$, which is absurd, since $d_p$ is convex. The assertion follows, and this proves $(i)$. If $K$ is unbounded, then we may show that it is a cylinder and therefore satisfies the Local Geodesic Property at every point, which proves $(ii)$ and completes the proof.\qed
\newhead{Existence - Review of Ye's Work}
\newsubhead{Overview}
\noindent Let $M:=M^{n+1}$ be a compact $(n+1)$-dimensional Riemannian manifold. We review the proof of existence of embedded constant curvature hypersurfaces in $M$ for sufficently large values of the curvature. We recall that in \cite{Ye}, Ye proves the existence of embedded hyperspheres of constant mean curvature near critical points of the scalar curvature function provided that these critical points are non-degenerate. Indeed, for any point $p\in M$, denote by $S_p\subseteq T_pM$ the sphere of unit vectors in $T_pM$, and for any $\varphi\in C^\infty(S_p)$ and any $t\in[0,\infty[$, we denote by $i^{t,\varphi}:S_p\rightarrow M$ the immersion given by:
\headlabel{HeadExistence}
\subheadlabel{SubHeadExistence}
$$
i^{t,\varphi}(y) = \opExp(t(1+\varphi(y))y).
$$
\noindent We call the triplet $(t,p,\varphi)$ a {\bf centred sphere}. We call $t$ the {\bf scale} of $(t,p,\varphi)$, and we call $p$ the {\bf centre} of $(p,\varphi)$. Let $(p_t)_{t\in[0,\epsilon[}$ be a smooth family of points in $M$, and let $(\varphi_t)_{t\in[0,\epsilon[}$ be a smooth family of smooth functions such that, for all $t$, $(t,p_t,\varphi_t)$ is a centred sphere with centre at $p_t$. We say that this family {\bf converges} to $p$ if and only if:
\medskip
\myitem{(i)} $p_0=p$; and
\medskip
\myitem{(ii)} $\varphi = O(t^2)$.
\medskip
\noindent The part of Ye's Theorem most relevant to our current context may now be stated as follows:
\proclaim{Theorem \nextprocno, {\bf Ye (1991)}}
\noindent If $p\in M$ is a non-degenerate critical point of the scalar curvature function, then there exists a smooth family, $(p_t,\varphi_t)_{t\in[0,\epsilon[}$, of centred spheres converging to $p$ such that, for all $t$, $(p_t,\varphi_t)$ has constant mean curvature equal to $t^{-1}$.
\endproclaim
\proclabel{ThmYe}
\noindent The hypotheses of Ye's Theorem are equivalent to the condition that the scalar curvature function be a Morse function. As outlined in the introduction, this is too restrictive for our purposes since it is not clear that metrics whose curvature functions are Morse Functions are themselves generic. Pacard and Xu overcome this difficulty in \cite{PacardXu}, thus proving the existence of constant mean curvature hyperspheres in general manifolds for sufficiently high values of curvature, but they still leave open the problem of estimating the spectrum of the Jacobi operators, which is our main objective. The following mild generalisation of Ye's Theorem allows us to easily bypass these difficulties:
\proclaim{Theorem \procref{ThmYeGeneralised}, {\bf Existence - Long Version}}
\noindent Let $R$ be the scalar curvature function of $M$, and let $f\in C^\infty(M)$ be a smooth real valued function. Let $p\in M$ be any point where the Hessian of $R_f:=R + 2(n+3)/(n+1)f$ is non-degenerate. There exists a smooth family, $(p_t,\varphi_t)_{t\in[0,\epsilon[}$, of centred spheres converging to $p$ such that, for all $t$, and for all $y\in S_{p_t}$, the $K$-curvature if $i^{t,\varphi_t}$ at the point $i^{t,\varphi_t}(y)$ is equal to:
$$
t^{-1}\left(1 + t^2(f\circ i^{t,\varphi_t})(y) - t^3\frac{(n+1)}{2(n+3)}\langle \nabla R_f, \tau_{p,p_t}y\rangle\right),
$$
\noindent where $\tau_{p,p_t}$ is the parallel transport from $p_t$ to $p$ along the unique minimising geodesic joining these two points.
\medskip
\noindent In particular, if $p$ is a non-degenerate critical point of $R_f$, then, for all $t\in[0,\epsilon[$, the $K$-curvature of $(p_t,\varphi_t)$ is prescribed by the function $t^{-1}(1+t^2 f)$.
\endproclaim
\remark In the sequel, we will refer to the centred spheres constructed by Theorem \procref{ThmYeGeneralised} as {\bf Ye's Spheres}.
\medskip
\noindent Theorem \procref{ThmYeGeneralised} is a relatively trivial generalisation of Theorem \procref{ThmYe}. We nonetheless present a fairly detailled discussion of Ye's proof, since the asymptotic expansions obtained here will be of importance in the sequel.
\newsubhead{Forms Over Radial Spheres}
\noindent We calculate various geometric forms over radial spheres. Since we closely follow Ye's approach - albeit from a slightly different angle - we will not provide detailed proofs where they do not enlighten.
\medskip
\noindent We recall the framework within which we work throughout the sequel. Let $M:=M^{n+1}$ be an $(n+1)$-dimensional Riemannian manifold. Given $p\in M$, we define $\Phi:\Bbb{R}^{n+1}\times\Bbb{R}^{n+1}\rightarrow M$ as follows: let $\opExp$ be the exponential mapping of $M$. Furnish $\Bbb{R}^{n+1}$ with the canonical metric and let $A:\Bbb{R}^{n+1}\rightarrow T_pM$ be an isometry. For all $x\in\Bbb{R}^{n+1}$, let $[0,x]$ be the straight line in $\Bbb{R}^{n+1}$ joining $0$ to $x$ and let $\tau_x$ be the parallel transport map in $M$ along the geodesic $\opExp(A\cdot[0,x])$ from $p=\opExp(A\cdot 0)$ to $\opExp(A\cdot x)$. We now define $\Phi_x(y) := \Phi(x,y)$ by:
$$
\Phi_x(y) = \opExp((\tau_x\circ A)\cdot y).
$$
\noindent We view $\Phi$ as a smooth family of exponential coordinate charts for $M$ near $p$ parametrised by the first component. For all $x$, we identify $M$ (locally) with the fibre $\left\{x\right\}\times\Bbb{R}^{n+1}$. This fibre carries the Euclidean metric, which we denote by $g^0$, as well as the pull-back of the Riemannian metric over $M$, which we denote by $g^x$.
\medskip
\noindent Now consider $[0,\infty[\times\Bbb{R}^{n+1}\times\Bbb{R}^{n+1}$. For all $(t,x)$, we define the metric $g^{t,x}$ over the fibre $\left\{(t,x)\right\}\times\Bbb{R}^{n+1}$ by:
$$
g^{t,x}(y) = g^x(ty) = t^{-2}(\delta_t^*g^x)(y),
$$
\noindent where $\delta_t:\Bbb{R}^{n+1}\rightarrow\Bbb{R}^{n+1}$ is the dilatation given by:
$$
\delta_t\cdot y = ty.
$$
\noindent In the sequel, we work with respect to an orthonormal basis, $e_1,...,e_{n+1}$ of the fibre chosen such that $e_{n+1}=y$, where $y$ is the position vector. The summation convention will be used whenever indices are repeated, a semi-colon appearing as a subscript indicates covariant differentiation, and the subscript, $y$, represents contraction with the position vector, and not an index. Finally, here and throughout the sequel, we adopt the following conventions:
$$
R_{ijkl} = \langle\nabla_{\partial_i}\nabla_{\partial_j}\partial_k - \nabla_{\partial_j}\nabla_{\partial_j}\partial_k - \nabla_{[\partial_i\partial_j]}\partial_k,\partial_l\rangle.
$$
\noindent and:
$$
\opRic_{ij} = -\frac{1}{n}\sum_{k=1}^nR_{kikj}\qquad\&\qquad R = \frac{1}{(n+1)}\sum_{k=1}^n\opRic_{kk}.
$$
\noindent These conventions are chosen such that the Ricci and scalar curvatures of the unit sphere in Euclidean space are equal to $\delta_{ij}$ and $1$ respectively. Denoting by $R^x_{ijkl}$ the Riemann curvature tensor of $M$ at the point $x$, using Jacobi fields, we readily obtain:
\proclaim{Lemma \nextprocno}
\noindent For all $(t,x,y)$:
$$
g^{t,x}_{ij}(y) = \delta_{ij} + \frac{t^2}{3}R^x_{yiyj} + \frac{t^3}{6}R^x_{yiyj;y} + \frac{t^4}{20}R^x_{yiyj;yy} + \frac{2t^4}{45}R^x_{yiyp}R^x_{yjyp} + O(t^5).
$$
\endproclaim
\proclabel{LemmaExpansionOfMetric}
\noindent We denote by $R^{t,x}$ the Riemann curvature tensor of the metric $g^{t,x}$, and we denote by $W^{t,x}_{ij}$ the contraction $R^{t,x}_{yiyj}$. We obtain:
\proclaim{Lemma \nextprocno}
\noindent For all $(t,x,y)$:
$$
W_{ij}^{t,x}(y) = t^2R^x_{yiyj} + t^3R^x_{yiyj;y} + \frac{t^4}{2}R^x_{yiyj;yy} + \frac{t^4}{3}R^x_{yiyp}R^x_{yiyp} + O(t^5).
$$
\endproclaim
\proclabel{LemmaExpansionOfCurvature}
\noindent Let $S\subseteq\Bbb{R}^{n+1}$ be the unit sphere centred on the origin. For all $(t,x)$, we denote by $S^{t,x}$ the unit sphere in the fibre $\left\{(t,x)\right\}\times\Bbb{R}^{n+1}$ with respect to the metric $g^0$. $S^{t,x}$ is trivially identified with $S$. Moreover, by classical differential geometry, for sufficiently small $t$, the radial lines in $\left\{(t,x)\right\}\times\Bbb{R}^{n+1}$ are also unit speed minimising geodesics with respect to the metric $g^{t,x}$. Consequently, $S^{t,x}$ is also the geodesic sphere in the fibre of radius $1$ with respect to the metric $g^{t,x}$, and, moreover, the unit length radial vector coincides with the unit normal vector over $S^{t,x}$ with respect to the metric $g^{t,x}$.
\medskip
\noindent Let $II^{t,x}$ and $A^{t,x}$ be the second fundamental form and the shape operator respectively of $S^{t,x}$ with respect to $g^{t,x}$. It is more natural to work with $2$-forms (such as $II^{t,x}$) than with linear maps (such as $A^{t,x}$), and since it is only necessary to work with the conjugacy class of $A^{t,x}$, here and in the sequel, we will identify it with:
$$
A^{t,x} = [II]^{t,x},
$$
\noindent where, for any bilinear form, $B^{t,x}$:
$$
[B]^{t,x}_{ij} = (g^{t,x})^{-1/2}_{ip}B_{pq}(g^{t,x})^{-1/2}_{qj}.
$$
\noindent This is meaningful, since $g^{t,x}$ will always lie in a small neighbourhood of the identity where the reciprocal of the square root is smoothly defined.
\proclaim{Lemma \nextprocno}
\noindent For all $(t,x)$, and at every point $y\in S$:
$$
II^{t,x}_{ij}(y) = \delta_{ij} + \frac{2t^2}{3}R^x_{yiyj} + \frac{5t^3}{12}R^x_{yiyj;y} + \frac{3t^4}{20}R^x_{yiyj;yy} + \frac{6t^4}{45}R^x_{yiyp}R^x_{yjyp} + O(t^5)
$$
\endproclaim
\proclabel{LemmaExpansionOfSecondFF}
\proof If $X:=X^n$ is a smoothly immersed hypersurface in an $(n+1)$-dimensional Riemannian manifold, $Y:=Y^{n+1}$, if $\msf{N}$ is the unit normal vector field over $X$ compatible with the orientations of $X$ and $Y$ and if $f:X\rightarrow\Bbb{R}$ is a smooth function such that $f\msf{N}$ corresponds to an infinitesimal normal perturbation of $X$, then the first order variation of the pull back of the metric of $X$ under this perturbation is given by:
$$
(D_f g)(U,V) = 2fII(U,V),
$$
\noindent where $II$ is the second fundamental form of $X$. In the current context, this yields:
$$
\partial_r g^{t,x}_{ry}(rX,rX)|_{r=1} = 2II^{t,x}(rX,rX).
$$
\noindent The result now follows by substituting the formula obtained in Lemma \procref{LemmaExpansionOfMetric} into the above relation.\qed
\proclaim{Lemma \nextprocno}
\noindent For all $(t,x)$, and at every point $y\in S$:
$$
A^{t,x}_{ij}(y) = \delta_{ij} + \frac{t^2}{3}R^x_{yiyj} + \frac{t^3}{4}R^x_{yiyj;y} + \frac{t^4}{10}R^x_{yiyj;yy} - \frac{t^4}{45}R^x_{yiyp}R^x_{yjyp} + O(t^5).
$$
\endproclaim
\proclabel{LemmaExpansionOfShapeOperator}
\proof This follows from Lemma \procref{LemmaExpansionOfMetric} and Lemma \procref{LemmaExpansionOfSecondFF}, since $A^{t,x} = [II]^{t,x}$.\qed
\medskip
\noindent For all $(t,x)$, let $\Gamma^k_{ij}$ be the relative Christophel Symbol of the Levi-Civita covariant derivative of $g^{t,x}$ with respect to that of $g^0$ (we suppress $(t,x)$ in this case to avoid excessively cumbersome notation). We are only interested in the diagonal elements of this tensor up to order $3$ in $t$.
\proclaim{Lemma \nextprocno}
\noindent For all $(t,x,y)$:
$$
\Gamma^k_{ii}(y) = -\frac{2t^2}{3}R^x_{yiki} - \frac{t^3}{2}R^x_{yiki;y} + \frac{t^3}{12}R^x_{yiyi;k} + O(t^4).
$$
\endproclaim
\proclabel{LemmaExpansionOfGamma}
\proof This follows by the Kohzul formula and Lemma \procref{LemmaExpansionOfMetric}.\qed
\medskip
\noindent We are interested in the $K$-curvature. Let $\opSymm(n)$ be the space of symmetric $n$-dimen\-sional real matrices. We define $\hat{K}:\opSymm(n)\rightarrow\Bbb{R}$ by:
$$
\hat{K}(A) = K(\lambda_1,...,\lambda_n),
$$
\noindent where $\lambda_1,...,\lambda_n$ are the eigenvalues of $A$. In the sequel, we suppress the hat and write $K$ in place of $\hat{K}$, where no ambiguity arises. Since $K$ is homogeneous of order $1$, we readily obtain:
\proclaim{Proposition \nextprocno}
\noindent At the identity:
\medskip
\myitem{$(i)$} $DK_\opId = \frac{1}{n}\opId$; and
\medskip
\myitem{$(ii)$} $D^2K_\opId(\cdot,\opId) = 0$.
\endproclaim
\proclabel{PropPropertiesOfK}
\noindent We now determine the $K$-curvature of $S^{t,x}$ with respect to $g^{t,x}$ for all $(t,x)$:
\proclaim{Proposition \nextprocno}
\noindent For all $(t,x)$ and for all $y\in S$, the $K$-curvature of $S^{t,x}$ with respect to $g^{t,x}$ at the point $y$ is given by:
$$\matrix
K^{t,x}(y) \hfill&= 1 - \frac{t^2}{3}\opRic^x_{yy} - \frac{t^3}{4}\opRic^x_{yy;y} -\frac{t^4}{10}\opRic^x_{yy;yy}-\frac{t^4}{45n}R^x_{ypyq}R^x_{ypyq}
\hfill\cr
&\qquad  + \frac{t^4}{18}D^2K_{\opId}(R^x_{y\cdot y\cdot},R^x_{y\cdot y\cdot}) + O(t^5).\hfill\cr
\endmatrix$$
\endproclaim
\proclabel{PropExpansionOfCurvature}
\proof By definition, $\hat{K}^{t,x}=K(\hat{A}^{t,x})$, and the result then follows by Taylor's Theorem, Lemma \procref{LemmaExpansionOfShapeOperator} and Proposition \procref{PropPropertiesOfK}.\qed
\newsubhead{Forms Over Perturbed Spheres}
\noindent For $\varphi\in C^\infty(S)$, we define $S^{t,x,\varphi}$ to be the graph of $t^2\varphi$ over $S^{t,x}$. We aim to calculate the $K$-curvature of $S^{t,x,\varphi}$ up to order $4$ in $t$. We first determine the first derivative of the metric, the second fundemental form and the shape operator with respect to normal perturbations up to order $1$ in $t$:
\proclaim{Proposition \nextprocno}
\noindent For any $f\in C^\infty(S)$:
$$
(Dg\cdot f) = 2f\delta_{ij} + O(t^2).
$$
\endproclaim
\proclabel{PropDiffOfMetric}
\proof Up to order $2$ in $t$, $g^{t,x}$ coincides with $g^0$, and thus, up to order $2$ in $t$, the variation of the restriction of $g^{t,x}$ with respect to normal perturbations coincides with that of $g^0$. The result now follows by the classical differential geometry of the unit sphere in Euclidean space.\qed
\medskip
\noindent In like manner, we obtain:
\proclaim{Proposition \nextprocno}
\noindent For any $f\in C^\infty(S)$:
$$
(DII\cdot f)_{ij} = \delta_{ij}f - f_{ij} + O(t^2).
$$
\endproclaim
\proclabel{PropDiffOfSecondFF}
\noindent Combining these relations yields:
\proclaim{Proposition \nextprocno}
\noindent For any $f\in C^\infty(S)$:
$$
(DA\cdot f)_{ij} = (-\delta_{ij}f - f_{ij}) + O(t^2).
$$
\endproclaim
\proclabel{PropDiffOfShapeOperator}
\proof This follows from Propositions \procref{PropDiffOfMetric} and \procref{PropDiffOfSecondFF} since $A = [II]$.\qed
\medskip
\noindent We now determine the first and second derivatives of the curvature operator up to order $2$ and order $0$ respectively in $t$. We say that an operator $E$ is {\bf even} if and only if it sends even functions to even functions.
\proclaim{Proposition \nextprocno}
\noindent For any $f\in C^\infty(S)$:
$$
DK\cdot f = -\frac{1}{n}(n+\Delta^0)f + t^2E_1\cdot f + O(t^3),
$$
\noindent where $E$ is even.
\endproclaim
\proclabel{PropDiffOfCurvature}
\proof Define the involution $I:\Bbb{R}^{n+1}\rightarrow\Bbb{R}^{n+1}$ by $I(x)=-x$. By Lemma \procref{LemmaExpansionOfMetric}:
$$
I^*g^{t,x} = g^{t,x} + O(t^3).
$$
\noindent For $f\in C^\infty(S^{t,x})$, let $K^{t,x,f}$ be the curvature with respect to $g^{t,x}$ of the graph of $f$ at the point over $x$. If $f$ is even, then, since $g^{t,x}$ is even up to order $2$ in $t$, for all $s$:
$$\matrix
&K^{t,x,sf}\circ I \hfill&= K^{t,x,sf} + O(st^3)\hfill\cr
\Rightarrow\hfill&\partial_sK^{t,x,sf}\circ I|_{s=0}\hfill&= \partial_s K^{t,x,sf}|_{s=0} + O(t^3)\hfill\cr
\Rightarrow\hfill&(DK\cdot f)\circ I\hfill&= DK\cdot f + O(t^3).\hfill\cr
\endmatrix$$
\noindent We conclude that $DK$ is even up to order $2$ in $t$. The explicit expansion of $DK$ up to order $1$ in $t$ follows immediately from Propositions \procref{PropPropertiesOfK} and \procref{PropDiffOfShapeOperator}, and this completes the proof.\qed
\proclaim{Proposition \nextprocno}
\noindent For any $f\in C^\infty(S)$:
$$
D^2K(f,f) = E_2(f,f) + O(t),
$$
\noindent where $E_2$ is even.
\endproclaim
\proclabel{PropSecondDiffOfCurvature}
\proof Taking $t=0$, we see that the second derivative of the $K$-curvature up to order $0$ in $t$ coincides with the second derivative of the $K$-curvature of the unit sphere in Euclidean space. Trivially this is invariant under the action of $O(n)$, and, in particular, is even. This completes the proof.\qed
\medskip
\noindent This allows us to determine the shape operator of $S^{t,x,\varphi}$. We first work up to order $3$ in $t$. For $y\in S^{t,x}=S$, we denote by $A^{t,x,\varphi}$ the shape operator of $S^{t,x,\varphi}$ at the point over $y$. We obtain:
\proclaim{Proposition \nextprocno}
\noindent For all $(t,x,\varphi)$, and for all $y\in S$:
$$
A^{t,x,\varphi}(y) = \delta_{ij} + t^2(\frac{1}{3}R^x_{yiyj} - \delta_{ij}\varphi - \varphi_{ij}) + \frac{t^3}{4}R^x_{yiyj;y} + O(t^4).
$$
\endproclaim
\proclabel{PropExpansionOfPertShapeOperator}
\noindent For $y\in S^{t,x}=S$, we denote by $K^{t,x,\varphi}$ the $K$-curvature of $S^{t,x,\varphi}$ at the point over $y$. We obtain:
\proclaim{Proposition \nextprocno}
\noindent For all $(t,x,\varphi)$ and for all $y\in S$:
$$
K^{t,x,\varphi}(y) = 1 - t^2(\frac{1}{3}\opRic^x_{yy} + \frac{1}{n}(n+\Delta^0)\varphi) - \frac{t^3}{4}\opRic^x_{yy;y} + O(t^4).
$$
\endproclaim
\proclabel{PropExpansionOfPertCurvature}
\noindent Finally, we denote by $K_4^{t,x,\varphi}$ the fourth order coefficient of the asymptotic expansion for $K^{t,x,\varphi}$. We do not require an explicit formula for $K_4^{t,x,\varphi}$. The following result suffices:
\proclaim{Proposition \nextprocno}
\noindent Suppose that $\varphi = \varphi_0 + O(t)$. If $\varphi_0$ is an even function, then so is $K_4^{t,x,\varphi}$.
\endproclaim
\proclabel{PropFourthOrderCoeffOfCurvature}
\proof By Propositions \procref{PropExpansionOfCurvature}, \procref{PropDiffOfCurvature} and \procref{PropSecondDiffOfCurvature}:
$$
K_4^{t,x,\varphi} = -\frac{1}{10}\opRic^x_{yy;yy} - \frac{1}{45n}R^x_{ypyq}R^x_{ypyq} + \frac{1}{18}D^2K_\opId(R^x_{y\cdot y\cdot},R^x_{y\cdot y\cdot}) + E_1(\varphi) + \frac{1}{2}E_2(\varphi,\varphi).
$$
\noindent The result now follows.\qed
\medskip
\noindent Finally, let $\msf{N}^{t,x,\varphi}$ be the unit normal vector field over $S^{t,x,\varphi}$. The following result will be of use in the sequel:
\proclaim{Proposition \nextprocno}
\noindent For all $(t,x,\varphi)$ and for all $y\in S$:
$$
\msf{N}^{t,x,\varphi}(y) = y - t^2\nabla^S\varphi + O(t^3),
$$
\noindent where $\nabla^S$ is the gradiant operator over $S^{t,x}$.
\endproclaim
\proclabel{PropPertNormal}
\remark Observe, in particular, that:
$$
\nabla^S\varphi = \nabla\varphi - \langle\nabla\varphi,y\rangle y.
$$
\proof We define $i^{t,x,\varphi}$ by:
$$
i^{t,x,\varphi}(y) = (1 + t^2\varphi(y))y,
$$
\noindent and we readily show that for any $X$ tangent to $S^{t,x}$:
$$
g^{t,x}(i^{t,x,\varphi}_*X,\msf{N}^{t,x,\varphi}) = O(t^3).
$$
\noindent The result follows.\qed
\newsubhead{Existence}
\noindent Using the asymptotic series obtained in the preceeding two sections, we now proceed to the proof of Theorem \procref{ThmYeGeneralised}. We recall from Section \headref{HeadPreliminaries} that, for all $k$, $\Pi_k:L^2(S)\rightarrow H_k$ is the orthogonal projection onto the $k$'th order spherical harmonics.
\proclaim{Proposition \nextprocno}
\noindent For all $x$:
$$
\Pi_1(\opRic^x_{yy;y}) = \frac{2(n+1)}{(n+3)}(\nabla R)^x_y.
$$
\endproclaim
\proclabel{PropProjectionOfCurvature}
\proof For all $p$, by Lemma \procref{LemmaFeynmanIntegrals}:
$$\matrix
\int_{S^n}\opRic^x_{yy;y}x_p\opdVol \hfill&= \opRic^x_{ij;k}\int_{S^n}x_ix_jx_kx_p\opdVol\hfill\cr
&= \frac{\Omega_n}{(n+1)(n+3)}((n+1)(\nabla R)^x_p + 2R_{ip;i}).\hfill\cr
\endmatrix$$
\noindent Applying the Second Bianchi Identity therefore yields:
$$
\int_{S^n}\opRic^x_{yy;y}x_p\opdVol = \frac{2\Omega_n}{(n+3)}(\nabla R)^x_p.
$$
\noindent However, by Lemma \procref{LemmaFeynmanIntegrals}:
$$
\int_{S^n}x_px_q\opdVol = \frac{\Omega_n}{(n+1)}\delta_{pq}.
$$
\noindent The result follows.\qed
\medskip
\noindent Define $i^{t,x,\varphi}:S\rightarrow M$ by:
$$
i^{t,x,\varphi} = \Phi_x(t(1+t^2\varphi(y))y).
$$
\noindent The mapping $i^{t,x,\varphi}$ is a parametrisation of $S^{t,x,\varphi}$, now viewed as a submanifold of $M$. Observe that the $K$-curvature of this immersed submanifold is equal to $t^{-1}K^{t,x,\varphi}$. For $f\in C^\infty(M)$, we consider the mapping $F:M\times C^\infty(S)\times [0,\infty[\rightarrow C^\infty(S)$ given by:
$$
F^{t,x,\varphi}(y) = 1 + t^2(f\circ i^{t,x,\varphi})(y).
$$
\noindent Taylor's Theorem readily yields:
\proclaim{Proposition \nextprocno}
\noindent For all $(t,x,y)$:
$$
F^{t,x,\varphi}(y) = 1 + t^2f^x + t^3 f^x_{;y} + \frac{1}{2}t^4f^x_{;yy} + O(t^5).
$$
\endproclaim
\proclabel{PropExpansionOfFunction}
\noindent Following \cite{Ye}, we obtain:
\proclaim{Proposition \nextprocno}
\noindent Let $f\in C^\infty(M)$ be smooth. For all $x$, there exist unique functions $\varphi_{0,x},\varphi_{1,x}\in H_1^\perp$ such that:
$$\matrix
\frac{1}{n}(n+\Delta^0)\varphi_{0,x} \hfill&= -\frac{1}{3}\opRic^x_{yy} - f^x,\hfill\cr
\frac{1}{n}(n+\Delta^0)\varphi_{1,x} \hfill&= -\frac{1}{4}\opRic^x_{yy;y} + \frac{(n+1)}{2(n+3)}R^x_{;y}.\hfill\cr
\endmatrix$$
\noindent Moreover, $\varphi_{0,x}$ and $\varphi_{1,x}$ are even and odd respectively.
\endproclaim
\proclabel{PropFirstAndSecondOrderTerms}
\proof Since the right hand side of the first line is an even function, it is an element of $H_1^\perp$. Likewise, by Proposition \procref{PropProjectionOfCurvature}, the right hand side of the second line is also an element of $H_1^\perp$. However, the restriction of $(n+\Delta^0)$ to $H_1^\perp$ is invertible, and the first result follows. Since $(n+\Delta^0)$ preserves the spaces of even and odd functions, so does its inverse, and it follows that $\varphi_{0,x}$ and $\varphi_{1,x}$ are even and odd respectively. This completes the proof.\qed
\medskip
\noindent We thus consider the mapping $\Psi:M\times H_1^\perp\times[0,\infty[\rightarrow C^\infty(S)$ given by:
$$
\Psi(t,x,\varphi) = \frac{1}{t^3}(K^{t,x,\hat{\varphi}} - F^{t,x,\hat{\varphi}}),
$$
\noindent where:
$$
\hat{\varphi} = \varphi_{0,x} + t(\varphi_{1,x} + \varphi).
$$
\noindent The function $t^3\Psi$ is trivially smooth. Observe, moreover, that, by Propositions \procref{PropExpansionOfPertCurvature} and \procref{PropExpansionOfFunction}, and by definition of $\varphi_{0,x}$ and $\varphi_{1,x}$, $t^3\Psi = O(t^3)$, and we deduce that $\Psi$ extends to a smooth function at $t=0$ such that:
$$
\Psi(0,0,0) = -\frac{(n+1)}{2(n+3)}R^0_{;y} - f^0_{;y}.
$$
\noindent Let $D_t\Psi$, $D_x\Psi$ and $D_\varphi\Psi$ be the partial derivatives of $\Psi$ with respect to the first second and third variables respectively. We readily obtain:
\proclaim{Proposition \nextprocno}
\noindent At $(0,0,0)$:
$$\matrix
D_x\Psi\cdot V \hfill&= -\frac{(n+1)}{2(n+3)}R^0_{;yV} - f_{;yV},\hfill\cr
D_\varphi\Psi\cdot g \hfill&= -\frac{1}{n}(n+\Delta^0)g,\hfill\cr
D_t\Psi\cdot \partial_t&= \psi,\hfill\cr
\endmatrix$$
\noindent where $\psi$ is an even function.
\endproclaim
\proclabel{PropDerivativeOfFunctional}
\proof The first formulae are trivial. The third formula follows from Propositions \procref{PropExpansionOfPertCurvature}, \procref{PropFourthOrderCoeffOfCurvature} and \procref{PropExpansionOfFunction}, bearing in mind that $\varphi_{0,x}$ is even.\qed
\proclaim{Proposition \nextprocno}
\noindent If the Hessian of $R_f:=R + 2(n+3)/(n+1)f$ is non-degenerate, then there exists a unique smooth family $(\varphi_t,x_t)$ such that $(\varphi_0,x_0)=(0,0)$, and:
$$
\Psi(t,x_t,\varphi_t) = -\frac{(n+1)}{2(n+3)}R^0_{;y} - f^0_{;y}.
$$
\noindent Moreover, $x(t)=O(t^2)$.
\endproclaim
\proclabel{PropYeGeneralised}
\proof By Proposition \procref{PropDerivativeOfFunctional}, at $(0,0,0)$, $D_x\Psi\oplus D_\varphi\Psi:T_pM\oplus H_1^\perp\rightarrow L^2(S)$ is a Banach-Space isomorphism. The result now follows by the Implicit Function Theorem for Banach spaces. Define $\Gamma\subseteq M\times H_1^\perp\times[0,\infty[$ by:
$$
\Gamma = \left\{(x,\varphi,t)\ |\ \Psi(x,\varphi,t) = -\frac{(n+1)}{2(n+3)}R^0_{;y} - f^0_{;y}\right\}.
$$
\noindent Then $\Gamma$ is a smooth $1$-dimensional submanifold, and $t\mapsto (x_t,\varphi_t,t)$ is a parametrisation of $\Gamma$ near $(0,0,0)$. Since $D_t\Psi$ is an even function, it is an element of $H_1^\perp$, and so there exists $f\in H_1^\perp$ such that:
$$
D\Psi\cdot(0,f,\partial_t) = 0.
$$
\noindent The vector $(0,f,\partial_t)$ is therefore tangent to $\Gamma$ at $(0,0,0)$, and since its first component vanishes, $x'(0)=0$, and the second assertion now follows.\qed
\medskip
\noindent Theorem \procref{ThmYeGeneralised} now follows:
\medskip
{\bf\noindent Proof of Theorem \procref{ThmYeGeneralised}:\ }This follows immediately from Proposition \procref{PropYeGeneralised}.\qed
\medskip
\noindent In addition, the uniqueness part of the Implicit Function Theorem for Banach Spaces immediately yields:
\proclaim{Corollary \nextprocno}
\noindent There exists $k,\alpha,\epsilon>0$ such that, if $\Psi(t,x,\varphi)=0$, $d(x,p)<\epsilon$, $t<\epsilon$ and:
$$
\|(t^2\varphi_{0,x} + t^3\varphi_{1,x}) - \varphi\|_{C^{k,\alpha}} \leqslant \epsilon^4,
$$
\noindent then $\varphi$ is one of Ye's spheres.
\endproclaim
\proclabel{CorUniqueness}
\newhead{The Signature of The Hypersurfaces}
\newsubhead{Overview}
\noindent We determine the signature of the hypersurfaces constructed in Section \headref{HeadExistence}. Let $K$ be a curvature function, let $f\in C^\infty(M)$ be a smooth function, and let $\Sigma\subseteq M$ be a compact immersed hypersurface of prescribed $K$-curvature equal to $f$. Let $JK$ be the Jacobi Operator of $K$ over $\Sigma$. The Jacobi Operator of the pair $(K,f)$ over $\Sigma$ is defined by:
$$
J(K,f):C^\infty(\Sigma)\rightarrow C^\infty(\Sigma);\psi\mapsto JK\cdot\psi - \langle\msf{N},\nabla f\rangle\psi,
$$
\noindent where $\msf{N}$ is the unit normal vector field over $\Sigma$ compatible with the orientation of $\Sigma$. As discussed in \cite{SmiDTI}, the spectrum of $J(K,f)$ is discrete and its real eigenvalues are bounded below. We say that $\Sigma$ is {\bf non-degenerate} if and only if $J(K,f)$ is non-degenerate, and, in this case, we define the {\bf signature} of $\Sigma$, which we denote by $\opSgn(\Sigma)$, by:
$$
\opSgn(\Sigma) = \opSgn(J(K,f)) = (-1)^{\opMult(J(K,f))},
$$
\noindent where $\opMult(J(K,f))$ is the number of strictly negative real eigenvalues of $J(K,f)$ counted with multiplicity. We show:
\proclaim{Theorem \procref{ThmSignature}, {\bf Signature}}
\noindent Let $M:=M^{n+1}$ be an $(n+1)$-dimensional Riemannian manifold. Let $R$ be the scalar curvature function of $M$ and let $f\in C^\infty(M)$ be a smooth function such that $R_f:=R + 2(n+3)/(n+1)f$ is a Morse Function. Let $p\in M$ be a critical point of $R_f$ and let $(p_t,\varphi_t)_{t\in[0,\epsilon[}$ be the family of Ye's spheres about $p$. Then, for sufficiently small $t$, $\Sigma_t:=(t,p_t,\varphi_t)$ is non-degenerate, and its signature is given by:
$$
\opSgn(\Sigma_t) = -(-1)^{n+1}\opSgn(\opHess(R_f)(p)).
$$
\endproclaim
\newsubhead{Perturbation Theory of Eigenvalues}
\noindent We consider a family, $(\Sigma_t)_{t\in[0,\epsilon[}=(p_t,\varphi_t)_{t\in[0,\epsilon[}$ of Ye's spheres about a critical point, $p$, of the function $R_f$. We first review the basic perturbation theory we will require to determine the signatures of the $(\Sigma_t)_{t\in[0,\epsilon[}$. For all $t$, we define the operator $J_t$ on $C^\infty(\Sigma_t)$ by:
\subheadlabel{SubHeadPerturbationTheory}
$$
J_t = t^2J(K,t^{-1}(1 + t^2f)).
$$
\noindent $J_t$ is $t^2$ times the Jacobi operator of $(K,t^{-1}(1+t^2f))$ over $\Sigma_t$. Bearing in mind the notation of Section \headref{HeadExistence}. We define $F^{t,p_t}$ by:
$$
F^{t,p_t}(y) = 1 + t^2f\circ\Phi(p_t,ty).
$$
\noindent Then, for all $t$, $J_t$ is also the Jacobi operator of the $K$-curvature of $(S^{t,p_t,\varphi_t},F^{t,p_t})$ with respect to the metric $g^{t,p_t}$. In particular, it is a smooth function of $t$. In the sequel, we need to consider the asymptotic expansion of $J_t$ up to order $4$ in $t$. We define the operators $A_0, ...,A_4$ as follows:
$$
J_t = A_0 + A_1t + A_2t^2 + A_3t^3 + A_4t^4 + O(t^5).
$$
\noindent We will see presently that:
\medskip
\myitem{(i)} $A_0$ is self adjoint, $\opKer(A_0)=H_1$, and the $A_0$ defines an isomorphism of $H_1^\perp$ onto itself;
\medskip
\myitem{(ii)} $A_1=0$;
\medskip
\myitem{(iii)} $A_2\cdot H_1\subseteq H_1^\perp$; and
\medskip
\myitem{(iv)} $A_3\cdot H_1\subseteq H_1^\perp$.
\medskip
\noindent This information yields a formula for the perturbation up to order $4$ in $t$ of the restriction of $J_t$ to the null eigenspace $H_1$. As in \cite{Kato} (but see also \cite{SmiDTI}), there exist smooth families, $(E_t)_{t\in[0,\epsilon[}$ and $(F_t)_{t\in[0,\epsilon[}$ of subspaces of $L^2(S)$ such that:
\medskip
\myitem{(i)} $E_0=H_1$ and $F_0=H_1^\perp$; and
\medskip
\myitem{(ii)} for all $t$, $J_t$ preserves $E_t$ and $F_t$.
\medskip
\noindent In particular, since $(E_t)_{t\in[0,\epsilon[}$ and $(F_t)_{t\in[0,\infty[}$ are continuous, we may assume that for all $t$:
$$
L^2(S) = E_t\oplus F_t.
$$
\noindent The eigenvalues of the restriction of $J_t$ to $E_t$ are thus the perturbations of the degenerate eigenvalue $l_1=0$ of $J_0 = A_0 = -\frac{1}{n}(\Delta^0 + n)$ (c.f. \cite{Kato}). Recall that $H_1$ is $(n+1)$-dimensional. Let $(v_t^1,...,v_t^{n+1})_{t\in[0,\epsilon[}$ be a smooth family such that, for all $t$, $(v_t^1,...,v_t^{n+1})$ constitutes an orthonormal basis for $E_t$ with respect to the $L^2$ metric. For all $1\leqslant i\leqslant (n+1)$, we denote:
$$
v_t^i = v_0^i + tv_1^i + t^2v_2^i + t^3v_3^i + t^4v_4^i + O(t^5).
$$
\noindent We may assume that $v_1^i\in H_1^\perp$ for all $1\leqslant i\leqslant (n+1)$. For all $1\leqslant i,j\leqslant (n+1)$, denote:
$$
a^{ij}_t = \langle J_tv_t^i,v_t^j\rangle.
$$
\noindent We denote:
$$
a^{ij}_t = a^{ij}_0 + ta^{ij}_1 + t^2a^{ij}_2 + t^3a^{ij}_3 + t^4a^{ij}_4 + O(t^5).
$$
\proclaim{Proposition \nextprocno}
\noindent For all $1\leqslant i,j\leqslant (n+1)$:
$$
a^{ij}_0 = a^{ij}_1 = a^{ij}_2 = a^{ij}_3 = 0,
$$
\noindent and:
$$
a^{ij}_4 = \langle A_4v_0^i,v_0^j\rangle + \langle A_2v_2^i,v_0^j\rangle,
$$
\noindent where, for all $i$:
$$
A_0v_2^i = A_2v_0^i.
$$
\endproclaim
\proclabel{PropExpansionOfRestriction}
\remark Since $A_2\cdot H_1\subseteq H_1^\perp$, we only need to know $v_2^i$ modulo $H_1$.
\medskip
\proof We begin by determining $v_t^i$ up to order $2$ in $t$. Let $(u_t)_{t\in[0,\epsilon[}$ be a smooth family of elements of $F_t$. Since $J_t$ preserves $E_t$ and $F_t$, for all $i,t$:
$$
\langle u_t, J_tv_t^i\rangle = 0.
$$
\noindent Denote $u_t=u_0 + tu_1 + t^2u_2 + O(t^3)$. Expanding the above relation yields, for all $k\geqslant 0$:
$$
\sum_{a+b+c=k}\langle u_a, A_bv_c^i\rangle = 0.
$$
\noindent Since $A_1=0$, and since the restriction of $A_0$ to $E_0=H_1$ vanishes, taking $k=1$ yields:
$$
u_0(A_0v_1^i)=0.
$$
\noindent Thus, since $u_0\in E_0^\perp$ is arbitrary:
$$
A_0v_1^i \in E_0.
$$
\noindent However, by hypotheses, $\opIm(A_0)\subseteq E_0^\perp$, and so $A_0v_1^i\in E_0^\perp\minter E_0$. Thus $A_0v_1^i=0$. Moreover, by hypothesis, $v_1^i\in E_0^\perp$, and so, since the restriction of $A_0$ to $E_0^\perp$ is injective, $v_1^i=0$ for all $i$. Taking $k=2$ now yields:
$$
u_2(A_0v_0^i) + u_0(A_2v_0^i) + u_0(A_0v_2^i) = 0.
$$
\noindent As before, $A_0v_0^i=0$, and so:
$$
u_0(A_2v_0^i + A_0v_2^i) = 0.
$$
\noindent Since $A_2v_0^i\in E_1^\perp$, and $A_0v_2^i\in\opIm(A_0)=E_1^\perp$, and since $u_0\in E_1^\perp$ is arbitrary:
$$
A_2v_0^i + A_0v_2^i = 0.
$$
\noindent This yields a satisfactory expansion for $v$ up to order $2$ in $t$.
\medskip
\noindent Expanding the relation $\langle J_tv_t^i,v_t^j\rangle = a^{ij}_t$ yields, for all $k$:
$$
a_k^{ij} = \sum_{(a,b,c)\in I_k^0}\langle A_av_b^i,v_c^j\rangle,
$$
\noindent where, for all $k$:
$$
I_k^0 = \left\{(a,b,c)\ |\ a,b,c\geqslant 0, a+b+c=k\right\}.
$$
\noindent Let $I_k$ be the set of all triplets in $(a,b,c)\in I_k^0$ such that:
$$
\langle A_av_b^i,v_c^j\rangle \neq 0.
$$
\noindent Trivially, for all $1\leqslant i,j\leqslant (n+1)$ and for all $k$:
$$
a_k^{ij} = \sum_{(a,b,c)\in I_k}\langle A_av_b^i,v_c^j\rangle,
$$
\noindent Since $A_1=0$ and $v_1^i=0$ for all $i$, $I_k\subseteq I_k^1$, where:
$$
I_k^1 = \left\{(a,b,c)\in I_k^0\ |\ a,b,c\neq 1\right\}.
$$
\noindent Since $A_0$ is self-adjoint and since its restriction to $E_0$ vanishes, $I_k\subseteq I_k^2$, where:
$$
I_k^2 = \left\{(a,b,c)\in I_k^1\ |\ a+b>0\ \&\ a+c>0\right\}.
$$
\noindent Finally, since $A_2E_0,A_3E_0\in E_0^\perp$, $I_k\subseteq I_k^3$, where:
$$
I_k^3 = I_k^2\setminus\left\{(2,0,0),(3,0,0)\right\}.
$$
\noindent Thus:
$$
I_0 = I_1 = I_2 = I_3 = \emptyset,
$$
\noindent and:
$$
I_4 = \left\{(4,0,0), (2,2,0), (2,0,2), (0,2,2)\right\}.
$$
\noindent Thus, for all $1\leqslant i,j\leqslant (n+1)$:
$$
a^{ij}_0 = a^{ij}_1 = a^{ij}_2 = a^{ij}_3 = 0,
$$
\noindent and:
$$
a_4^{ij} = \langle A_4v_0^i,v_0^j\rangle + \langle A_2v_2^i,v_0^j\rangle + \langle A_2v_0^i,v_2^j\rangle + \langle A_0v_2^i,v_2^j\rangle.
$$
\noindent Since $A_2v_0^i + A_0v_2^j=0$:
$$
a_4^{ij} = \langle A_4v_0^i,v_0^j\rangle + \langle A_2v_2^i,v_0^j\rangle.
$$
\noindent This completes the proof.\qed
\newsubhead{Coefficients up to Order $2$}
\noindent We recall that, for all $t$, $J_t$ is the Jacobi Operator of $(S^{t,p_t,\varphi_t},F^{t,p_t})$ with respect to the metric $g^{t,p_t}$. This allows us to readily determine $A_0$ and $A_1$:
\proclaim{Proposition \nextprocno}
\myitem{(i)} $A_0 = -\frac{1}{n}(\Delta_0 + n)$; and
\medskip
\myitem{(ii)} $A_1 = 0$.
\endproclaim
\proclabel{PropZeroethAndFirstOrderCoefficients}
\proof Up to order $1$ in $t$, $g^{t,x}$ coincides with the Euclidean metric $g^0$. Likewise, up to order $1$ in $t$, the sphere $S^{t,x,\varphi}$ coincides with the sphere of unit radius. Finally, up to order $2$ in $t$, $F$ is constant. Thus, up to order $1$ in $t$:
$$
J_t = -\frac{1}{n}(\Delta_0 + n) + O(t^2).
$$
\noindent The result now follows.\qed
\medskip
\noindent In order to determine the second order term, it is useful to explicitly determine $\varphi$ up to order $2$ in $t$. Bearing in mind the notation of Proposition \procref{PropFirstAndSecondOrderTerms}:
\proclaim{Proposition \nextprocno}
\noindent For all $x$ and for all $y\in S$:
$$
\varphi_{0,x}(y) = \frac{n}{3(n+2)}\opRic^x_{yy} - \frac{2(n+1)}{3(n+2)}R^x - f^x.
$$
\endproclaim
\proclabel{PropSecondOrderSoln}
\proof Trivially:
$$
-\frac{1}{3}\opRic^x_{yy}-f^x = -\frac{1}{3}(\opRic^x_{yy} - R^x) - (\frac{1}{3}R^x + f^x),
$$
\noindent where:
$$
(\opRic^x_{yy} - R^x) \in H_2,\qquad \frac{1}{3}R^x + f^x\in H_0.
$$
\noindent Since $H_0$ and $H_2$ are eigenspaces of $\frac{1}{n}(\Delta^0 + n)$, the result now follows by dividing each term by its corresponding eigenvalue.\qed
\medskip
\noindent We will see that it suffices to know the action of $A_2$ on $H_1$:
\proclaim{Proposition \nextprocno}
\noindent If $\psi\in H_1$, then:
$$
A_2\psi = -\frac{4(n+3)}{3(n+2)}\Pi_3(\opRic_{yy}^x\psi).
$$
\endproclaim
\proclabel{PropSecondOrderOperator}
\proof By Lemma \procref{LemmaExpansionOfMetric}:
$$
g^{t,x,\varphi}_{ij} = \delta_{ij} + t^2(\frac{1}{3}R^x_{yiyj} + 2\varphi\delta_{ij}) + O(t^3).
$$
\noindent By Lemma \procref{LemmaExpansionOfCurvature}:
$$
W^{t,x,\varphi}_{ij} = t^2R^x_{yiyj} + O(t^3).
$$
\noindent By Lemma \procref{LemmaExpansionOfSecondFF} and Proposition \procref{PropDiffOfSecondFF}:
$$
II^{t,x,\varphi}_{ij} = \delta_{ij} + t^2(\frac{2}{3}R^x_{yiyj} + \varphi\delta_{ij} - \varphi_{ij}) + O(t^3).
$$
\noindent Thus:
$$
[W]^{t,x,\varphi}_{ij} - (A^{t,x,\varphi})^2_{ij} = -\delta_{ij} + t^2(\frac{1}{3}R^x_{yiyj} + 2\varphi\delta_{ij} + 2\varphi_{ij}) + O(t^3).
$$
\noindent Thus, by Propositions \procref{PropPropertiesOfK} and \procref{PropFirstAndSecondOrderTerms}:
$$\matrix
DK_{A^{t,x,\varphi}}([W]^{t,x,\varphi} - (A^{t,x,\varphi})^2)\psi \hfill&=
-\psi + t^2(-\frac{1}{3}\opRic^x_{yy} + \frac{2}{n}(n+\Delta^0)\varphi)\psi + O(t^3)\hfill\cr
&=-\psi - t^2\opRic^x_{yy} - 2t^2f^x\psi + O(t^3).\hfill\cr
\endmatrix$$
\noindent We now calculate the Hessian of $\psi$. We extend $\psi\in H_1$ to an homogeneous function of order $1$, $\hat{\psi}\in\hat{\Cal{H}}_1$. Observe that the restriction of $\hat{\psi}$ to $S^{t,x,\varphi}$ at the point over $y\in S^{t,x}$ is equal to $(1+t^2\varphi(y))\psi(y)$. Let $\msf{N}^{t,x,\varphi}$ be the outward pointing unit normal over $S^{t,x,\varphi}$ with respect to the metric $g^{t,x,\varphi}$. By Proposition \procref{PropPertNormal}, expressing $\msf{N}^{t,x,\varphi}$ with respect to its tangential and normal components, we obtain:
$$
\msf{N}^{t,x,\varphi}= (-t^2\nabla^S\varphi,1) + O(t^3),
$$
\noindent where $\nabla^S$ is the gradient operator of $S$. Thus:
$$
\langle\nabla\hat{\psi},\msf{N}^{t,x,\varphi}\rangle = \psi - t^2\langle\nabla^S\varphi,\nabla^S\hat{\psi}\rangle + O(t^3).
$$
\noindent Let $\opHess^0$ and $\opHess^{t,x}$ be the Hessian operators of $g^0$ and $g^{t,x}$ respectively. Trivially:
$$
\opHess^0(\hat{\psi})_{ii} = 0.
$$
\noindent By Lemma \procref{LemmaExpansionOfGamma}, the diagonal elements of $\opHess^{t,x}(\hat{\psi})$ are therefore:
$$
\opHess^x(\hat{\psi})_{ii} = \frac{2t^2}{3}R^x_{yiki}\psi_k + O(t^3).
$$
\noindent Thus, by Lemma \procref{LemmaHessianOfRestrictionGeneral}, the diagonal elements of the Hessian of the restriction of $\hat{\psi}$ to $S^{t,x,\varphi}$ are:
$$
\hat{\psi}_{ii} = \frac{2t^2}{3}R^x_{yiki}\psi_k - (\psi - t^2\langle\nabla^S\varphi,\nabla^S\psi\rangle)II^{t,x,\varphi}_{ii} + O(t^3).
$$
\noindent By definition, up to order $2$ in $t$:
$$
DK_{A^{t,x,\varphi}}(A^{t,x,\varphi}) = K(A^{t,x,\varphi}) = 1 + t^2f^x + O(t^3).
$$
\noindent Thus, bearing in mind that, over $S^{t,x,\varphi}$, $\hat{\psi} = (1+t^2\varphi)\psi$, this yields:
$$\matrix
(1+t^2\varphi)DK_{A^{t,x,\varphi}}([\opHess(\psi)]_{ii}) \hfill&= -(\psi - t^2\langle\nabla^S\varphi,\nabla^S\psi\rangle)DK_{A^{t,x,\varphi}}(A^{t,x,\varphi}_{ii})\hfill\cr
&\qquad + t^2(-\frac{2}{3}\opRic^x_{yk}\psi_k - \frac{1}{n}(\Delta^S\varphi)\psi - \frac{2}{n}\varphi_i\psi_i) + O(t^3)\hfill\cr
&=-\psi - t^2(\frac{2}{3}\opRic^x_{yk}\psi_k + f^x\psi + \frac{1}{n}(\Delta^S\varphi)\psi \hfill\cr
&\qquad - \frac{(n-2)}{n}\langle\nabla^S\varphi,\nabla^S\psi\rangle) + O(t^3).\hfill\cr
\endmatrix$$
\noindent Thus, bearing in mind Proposition \procref{PropSecondOrderSoln} and the fact that $f$ is constant up to order $2$ in $t$:
$$\matrix
J_t\cdot\psi \hfill&= t^2(-\frac{4}{3}\opRic^x_{yy}\psi + \frac{2}{3}\opRic^x_{yk}\psi_k - 2(\varphi + f^x)\psi - \frac{(n-2)}{n}\langle\nabla^S\varphi,\nabla^S\psi\rangle) + O(t^3)\hfill\cr
&=t^2( -\frac{4(n+3)}{3(n+2)}\opRic^x_{yy}\psi + \frac{8}{3(n+2)}\opRic^x_{yk}\psi + \frac{4(n+1)}{3(n+2)}R^x\psi) + O(t^3).\hfill\cr
\endmatrix$$
\noindent Observe that the second order term is the restriction to $S$ of a polynomial function of order $3$. Thus, by Lemma \procref{LemmaHarmonicDecomposition}:
$$
J_t\psi + O(t^3) \in H_1\oplus H_3.
$$
\noindent However, by Lemma \procref{LemmaFeynmanIntegrals}:
$$\matrix
\int_{S}\opRic^x_{yy}\psi x_l\opdVol \hfill&= \opRic^x_{ij}\psi_k\int_{S}x_ix_jx_kx_l\opdVol\hfill\cr
&= \frac{\Omega_n}{(n+1)(n+3)}\opRic^x_{ij}\psi_k(\delta_{ij}\delta_{kl} + \delta_{ik}\delta_{jl} + \delta_{il}\delta_{jk})\hfill\cr
&= \frac{\Omega_n}{(n+3)}R^x\psi_l + \frac{2\Omega_n}{(n+1)(n+3)}\opRic_lk\psi_k.\hfill\cr
\endmatrix$$
\noindent Moreover, by Lemma \procref{LemmaFeynmanIntegrals} again:
$$
\int_{S}x_ix_j\opdVol = \frac{\Omega_n}{(n+1)}\delta_{ij}.
$$
\noindent Thus:
$$
\Pi_1(\opRic^x_{yy}\psi) = \frac{(n+1)}{(n+3)}R^x\psi + \frac{2}{(n+3)}\opRic_{yk}\psi_k.
$$
\noindent Consequently:
$$\matrix
J_t\cdot\psi\hfill&= -\frac{4(n+3)}{3(n+2)}\opRic^x_{yy}\psi + \frac{4(n+3)}{3(n+2)}\Pi_1(\opRic^x_{yy}\psi)\hfill\cr
&=-\frac{4(n+3)}{3(n+2)}\Pi_3(\opRic^x_{yy}\psi).\hfill\cr
\endmatrix$$
\noindent This completes the proof.\qed
\newsubhead{Coefficients up to Order $4$}
\noindent We use a straightforward geometric argument to bypass explicitly calculating the fourth order perturbation of the Jacobi Operator. We proceed as follows: let $p\in M$ be a critical point of $R_f = R + \frac{2(n+3)}{(n+1)}f$, let $V$ be a unit vector at $p$, let $s\mapsto p_s$ be the unit speed geodesic leaving $p$ in the direction of $V$ and let $(p_{s,t},\varphi_{s,t})_{t\in[0,\epsilon[}$ be the family of Ye's spheres constructed about the point $p_s$. For two points $p,q\in M$ sufficiently close, let $\tau_{q,p}$ be the parallel transport from $p$ to $q$ along the unique minimising geodesic joining $p$ to $q$. We denote:
$$
\Psi_{s,t}(y) = (\tau_{p_{s,t},p_s}\circ\tau_{p_s,p})(y).
$$
\noindent Throughout the rest of this section, we identify $\varphi_{s,t}$ with the composition $\varphi_{s,t}\circ\Psi_{s,t}$. For all $(s,t)$, define $i_{s,t}:S\rightarrow M$ by:
$$
i_{s,t} = \opExp_{p_{s,t}}(t(1+t^2\varphi_{s,t})(y)\cdot \Psi_{s,t}(y)),
$$
\noindent and for all $(s,t)$ and for all $y\in S$, let $K_{s,t}(y)$ be the $K$ curvature of $\Sigma_{s,t}:=(i_{s,t},S)$ at the point $i_{s,t}(y)$. We define:
$$
F_{s,t} = t^{-1}(1 + t^2(f\circ i_{s,t})(y)),
$$
\noindent and we define:
$$
\hat{K}_{s,t}(y) = t(K_{s,t}(y) - F_{s,t}(y)).
$$
\noindent By construction:
$$\matrix
\hat{K}_{s,t}(y) \hfill&= -t^3\frac{(n+1)}{2(n+3)}(\nabla R_f^{p_s}\circ\tau_{p_s,p})(y)\hfill\cr
&= -st^3\frac{(n+1)}{2(n+3)}\nabla^2R_f^p(y,V) + O(s^2t^3).\hfill\cr
\endmatrix$$
\noindent For all $s,t$ sufficiently small, we denote, for $y\in T_pM$:
$$
j_{s,t}(y) = i_{st,t}(y).
$$
\noindent We aim to calculate the infinitesimal variation of $j_{s,t}$ with respect to $s$ up to order $3$ in $t$. We denote this by $\xi_t$:
$$
\xi_t(y) = \partial_s j_{s,t}(y)|_{s=0}.
$$
\noindent We introduce some notation. For $p\in M$, $X\in T_pM$, we denote:
$$
p + X := \opExp_p(X).
$$
\noindent Likewise, for $X,Y\in T_pM$, we denote:
$$
(p+X)+Y := (p+X) + \tau_{p+X,p}(Y).
$$
\noindent Finally, for $p,q\in M$, and for any quantity, $C$, we denote:
$$
p = q\ \opMod\ C,
$$
\noindent if and only $d(p,q)=O(C)$.
\medskip
\noindent Since $(p+X)+Y$ and $p+(X+Y)$ are smooth functions of $p$, $X$ and $Y$ which coincide when one of $X$ or $Y$ vanishes:
$$
(p + X) + Y = p + (X + Y)\ \opMod\ O(\|X\|\|Y\|).
$$
\noindent In particular, denoting $q=p+X$, this immediately yields:
$$
q + \tau_{q,p}Y = p + Y\ \opMod\ d(p,q)
$$
\noindent Moreover, using Jacobi fields, we readily obtain:
$$
(p + X) + Y = (p + Y) + (X + \frac{1}{2}R_{YX}Y)\ \opMod\ O(\|X\|^2) + O(\|X\|\|Y\|^3).
$$
\noindent Now denote:
$$
j^1_{s,t}(y) = (\opExp_{p_{st,t}}\circ\Psi_{st,t})(ty),
$$
\noindent so that, in particular:
$$
j_{s,t}(y) = j^1_{s,t}((1+t^2\varphi_{s,t}(y))y).
$$
\proclaim{Proposition \nextprocno}
\noindent There exists a vector $V'\in T_pM$ such that, up to order $3$ in $t$:
$$
\partial_s j^1_{s,t}(y)|_{s=0} = (\tau_{p_t + ty,p_t}\tau_{p_t,p})(tV + t^3V' + \frac{1}{2}t^3R_{yV}y) + O(t^4).
$$
\endproclaim
\proclabel{PropNormalVariation}
\proof For all $s$ let $c_s:\Bbb{R}\rightarrow T_{p_s}M$ be such that, for all $t$:
$$
p_{s,t} = p_s + c_s(t).
$$
\noindent By construction (c.f. Theorem \procref{ThmYeGeneralised}), for all $s$:
$$
c_s(0) = c_s'(0) = 0.
$$
\noindent We define the vector field $U$ such that, for all $s$, $U_s = 2c_s''(0)$. Then, modulo $O(st^4)+O(s^2t^2)$:
$$\matrix
&c_{st}(t)\hfill&=\tau_{p_{st},p}(c_p(t) + st^3\nabla_{V}U)\hfill\cr
\Rightarrow\hfill&p_{st,t}\hfill&=p_{st} + \tau_{p_{st},p}(c_p(t) + st^3\nabla_{V}U)\hfill\cr
& &=(p + stV) + (c_p(t) + st^3\nabla_{V}U)\hfill\cr
& &=(p + (c_p(t) + st^3\nabla_{V}U)) + (stV + \frac{1}{2}stR_{c_p(t)V}c_p(t))\hfill\cr
& &=(p + (c_p(t) + st^3\nabla_V U)) + stV.\hfill\cr
\endmatrix$$
\noindent However, modulo $O(st^4) + O(s^2t^2)$:
$$\matrix
(p+(c_p(t) + st^3\nabla_V U)) \hfill&= (p+c_p(t)) + st^3\nabla_V U\hfill\cr
&= p_{0,t} + st^3\tau_{p_{0,t},p}\nabla_V U.\hfill\cr
\endmatrix$$
\noindent Moreover, consider the function $\Psi_{s,t}$ given by:
$$
\Psi_{s,t} = \tau_{p,p_{0,t}}\circ\tau_{p_{0,t},p^2_{s,t}}\circ\tau_{p^2_{s,t},p^1_{s,t}}\circ\tau_{p^1_{s,t},p},
$$
\noindent where:
$$\matrix
p^1_{s,t} \hfill= p + (c_p(t) + s\nabla_V U)\hfill\cr
p^2_{s,t} \hfill= (p + c_p(t)) + s\nabla_V U.\hfill\cr
\endmatrix$$
\noindent Then, modulo $O(st^4) + O(s^2t^2)$:
$$
p_{st,t} = (p_{0,t} + st^3\tau_{p_{0,t},p}\nabla_V U) + st(\tau_{p_{0,t},p}\circ\Phi_{st^3,t})V.
$$
\noindent However, since $\Phi_{s,t}$ is a linear function which equals the identity when either $s$ or $t$ vanishes, for any vector $U$:
$$
\Phi_{s,t}U = U + O(st\|U\|).
$$
\noindent Thus, modulo $O(st^4) + O(s^2t^2)$:
$$\matrix
p_{st,t} \hfill&= (p_{0,t} + st^3\tau_{p_{0,t},p}\nabla_VU) + st\tau_{p_{0,t},p}V\hfill\cr
&= p_{0,t} + \tau_{p_{0,t},p}(stV + st^3\nabla_{V}U).\hfill\cr
\endmatrix$$
\noindent Consider now $\Phi_{A,B,C}(X)$ given by:
$$
\Phi_{A,B,C}(X) = (\tau_{p,p+C}\circ\tau_{p+C,(p+A)+B}\circ\tau_{(p+A)+B,p+A}\circ\tau_{p+A,p})(X).
$$
\noindent This is a smooth function of $A$, $B$, $C$ and $X$ which is linear in $X$ and equal to $X$ when one $B=C$ and one of $A$ or $B=C$ vanishes. Thus:
$$
\Phi_{A,B,C}(X) = X + O(\|A\|\|B\|\|X\|) + O(\|A-B\|\|X\|).
$$
\noindent Composing with $\tau_{(p+A)+B,p+C}\circ \tau_{p+C,p}$ yields:
$$
(\tau_{(p+A)+B,p+A}\circ\tau_{p+A,p})X = (\tau_{(p+A)+B,p+C}\circ\tau_{p+C,p})X + O(\|U\|\|V\|\|X\|).
$$
\noindent Substituting $A=stV$, $B=\tau_{p,p+A}c_{ts}(t)$, $C=c_0(t)$ and $X=ty$ yields, modulo $O(st^4)+O(s^2t^2)$:
$$\matrix
p_{st,t} + (\tau_{p_{st,t},p_{st}}\circ\tau_{p_{st},p})ty \hfill&= p_{st,t} + (\tau_{p_{st,t},p_{0,t}}\circ\tau_{p_{0,t},p})ty\hfill\cr
&=(p_{0,t} + \tau_{p_{0,t},p}(st^3\nabla_V U + stV)) + \tau_{p_{0,t},p}ty\hfill\cr
&=(p_{0,t} + \tau_{p_{0,t},p}ty) + \tau_{p_{0,t},p}(stV + st^3\nabla_{V}U + \frac{1}{2}st^3R_{yV}y).\hfill\cr
\endmatrix$$
\noindent Differentiating with respect to $s$ therefore yields:
$$
\partial_s j^1_{s,t}(y)|_{s=0} = \tau_{p_t + ty,p_t}\tau_{p_t,p}(tV + t^3V' + \frac{1}{2}t^3R_{yV}y) + O(t^4).
$$
\noindent where $V'=\nabla_V U$. This completes the proof.\qed
\proclaim{Proposition \nextprocno}
\noindent Up to order $3$ in $t$:
$$
\xi_t(y) = \tau_{p_t + ty,p_t}\tau_{p_t,p}(tV + t^3V' + \frac{1}{2}t^3 R_{y V}y) + O(t^4).
$$
\endproclaim
\proclabel{PropNormalVariationOfYeSpheres}
\proof By definition:
$$
j_{s,t}(y) = j^1_{s,t}((1+t^2\varphi_{st,t}(y))y).
$$
\noindent Thus, by Proposition \procref{PropNormalVariation}, the chain rule and the product rule:
$$
\partial_s j_{s,t}(y)|_{s=0} = \tau_{p_t + ty,p_t}\tau_{p_t,p}(tV + t^3V' + \frac{1}{2}t^3 R_{y V}y + t^3\partial_s\varphi_{st,t}|_{s=0}y) + O(t^4).
$$
\noindent However:
$$
\partial_s\varphi_{st,t}|_{s=0} = O(t).
$$
\noindent The result now follows.\qed
\medskip
\noindent For all $t$, let $\msf{N}_t(y)$ be the outward pointing unit normal vector of $\Sigma_{0,t}$ at the point over $y$.
\proclaim{Proposition \nextprocno}
\noindent Up to order $3$ in $t$:
$$
\langle\msf{N}_t(y),\xi_t(y)\rangle - t^3\langle\nabla\varphi,y\rangle \langle y,V\rangle + O(t^4)\in H_1.
$$
\endproclaim
\proclabel{PropNormalComponentOfVariation}
\proof By Proposition \procref{PropPertNormal}:
$$
\msf{N}_t(y) = \tau_{p_t + ty,p_t}\circ\tau_{p_t,p}(y - t^2(\nabla\varphi - \langle\nabla\varphi,y\rangle y)) + O(t^3).
$$
\noindent Thus, by Proposition \procref{PropNormalVariationOfYeSpheres}:
$$\matrix
\langle\msf{N}_t(y),\partial_s j_{s,t}(y)|_{s=0}\rangle \hfill&= t\langle V,y\rangle + t^3\langle V', y\rangle + \frac{1}{2}t^3\langle R_{y V}y,y\rangle\hfill\cr
&\qquad - t^3\langle\nabla^0\varphi,V\rangle + t^3\langle\nabla^0\varphi,y\rangle\langle y,V\rangle + O(t^4).\hfill\cr
\endmatrix$$
\noindent However:
$$
\langle R_{y V}y,y\rangle = 0,
$$
\noindent and:
$$
\langle V,y\rangle,\langle V',y\rangle \in H_1.
$$
\noindent Moreover, since $\varphi_0\in\Cal{H}_0\oplus\Cal{H}_2$ along $S$:
$$
\langle\nabla\varphi_0,V\rangle\in H_1.
$$
\noindent The result follows.\qed
\medskip
\noindent This yields:
\proclaim{Proposition \nextprocno}
\noindent For all $\psi\in H_1$:
\medskip
\myitem{(i)} $A_2\cdot\psi, A_3\cdot\psi\in H_1^\perp$; and
\medskip
\myitem{(ii)}
$$
A_4\psi + 2A_2\psi + \frac{(n+1)}{2(n+3)}(\nabla^2R_f)(y,\nabla \psi) + O(t^5) \in H_1^\perp.
$$
\endproclaim
\proclabel{PropFourthOrderOperator}
\proof Let $V\in T_pM$ be such that $\psi(y)=\langle y, V\rangle$, and so, in particular, $\nabla\psi=V$. By construction:
$$
\hat{K}_{st,t}(y) = -st^4\frac{(n+1)}{2(n+3)}(\nabla^2 R_f)(y,V) + O(st^5).
$$
\noindent Since $y_s$ is parallel along $\gamma$, by definition of $J_t$ (see Section \subheadref{SubHeadPerturbationTheory}), differentiating with respect to $s$ yields:
$$
(t^{-2}J_t)\langle\msf{N}_t(y),\xi_t(y)\rangle + t^{-1}\xi_t^T(1+t^3(\nabla R_f)(y)) = -t^3\frac{(n+1)}{2(n+3)}(\nabla^2R_f)(y,V) + O(t^4),
$$
\noindent where $\xi_t^T$ is the tangential component of $\xi_t$, considered as a first order differential operator. However, $\nabla R_f$ vanishes at $s=0$, and so:
$$
J_t\langle\msf{N}_t(y),t^{-1}\xi_t(y)\rangle = -t^4\frac{(n+1)}{2(n+3)}(\nabla^2R_f)(y,V) + O(t^5)
$$
\noindent By definition $\opIm(A_0)=H_1^\perp$. Thus, since $\msf{N}_t(y) = y + O(t^2)$ and $\xi_t(y)=tV + O(t^3)$, we obtain:
$$\matrix
&J_t\psi \hfill&= O(t^4)\hfill\cr
\Rightarrow\hfill& A_0\psi + t^2A_2\psi + t^3A_3\psi + O(t^4)\hfill&\in H_1^\perp\hfill\cr
\Rightarrow\hfill& A_2\psi,A_3\psi \hfill&\in H_1^\perp.\hfill\cr
\endmatrix$$
\noindent This proves the first assertion. Moreover, up to order $4$ in $t$, bearing in mind Proposition \procref{PropNormalComponentOfVariation}:
$$
J_t\langle\msf{N}_t(y),\xi_t(y)\rangle - t^4(A_4\psi + A_2\langle y,\nabla\varphi\rangle\psi) + O(t^5)\in H_1^\perp.
$$
\noindent However, by the explicit formula for $\varphi$:
$$
\langle y,\nabla\varphi\rangle\psi - 2\varphi\psi \in H_1.
$$
\noindent The result now follows.\qed
\medskip
\noindent This allows us to conclude:
\proclaim{Proposition \nextprocno}
\noindent For all $i,j$:
$$
a^{ij}_4 = \frac{(n+1)}{2(n+3)}(\nabla^2R_f)_{ij}.
$$
\endproclaim
\proclabel{PropFourthOrderCoeffOfJacobiOperator}
\proof Choose $\psi$ such that $v_0^i=\nabla\psi$. By Propositions \procref{PropExpansionOfRestriction} and \procref{PropSecondOrderOperator}:
$$
A_0 v_2^i = A_2 v_0^i = -\frac{4(n+3)}{3(n+2)}\Pi_3(\opRic_{yy}^x\psi).
$$
\noindent Thus:
$$
v_2^i = \frac{2n}{3(n+2)}\Pi_3(\opRic_{yy}^x\psi) + \psi',
$$
\noindent for some $\psi'\in H_1$. Thus, by Proposition \procref{PropSecondOrderSoln}:
$$
v_2^i - 2\varphi_0\psi\in H_1.
$$
\noindent Consequently:
$$\matrix
&A_2v_2^i - 2A_2\varphi\psi\hfill&\in H_1^\perp\hfill\cr
\Rightarrow\hfill&\langle A_2v_2^i,v_0^j\rangle \hfill&= 2\langle A_2\varphi \psi,v_0^j\rangle.\hfill\cr
\endmatrix$$
\noindent Thus, by Propositions \procref{PropExpansionOfRestriction} and \procref{PropFourthOrderOperator}:
$$\matrix
a^{ij}_4\hfill&=\langle A_4v_0^i,v_0^j\rangle + \langle A_2v_2^i,v_0^j\rangle\hfill\cr
&=\langle A_4v_0^i + 2A_2\varphi\psi,v_0^j\rangle\hfill\cr
&=-\frac{(n+1)}{2(n+3)}\langle(\nabla^2 R_f)(y,\nabla\psi),v_0^j\rangle.\hfill\cr
\endmatrix$$
\noindent The result follows.\qed
\medskip
\noindent We now prove Theorem \procref{ThmSignature}:
\medskip
{\bf\noindent Proof of Theorem \procref{ThmSignature}:\ }By \cite{Kato}, there exists $\delta>0$ such that, for $t$ sufficiently small, the eigenvalues of $J_t$ which lie in $B_\delta(0)$ coincide, with multiplicity, with the eigenvalues of $(a^{ij}_t)$.
\medskip
\noindent By Proposition \procref{PropExpansionOfRestriction}:
$$
a_0^{ij} = a_1^{ij} = a_2^{ij} = a_3^{ij} = 0.
$$
\noindent Moreover, by Proposition \procref{PropFourthOrderCoeffOfJacobiOperator}, and by definition of $f$, $a_4^{ij}$ is non-degenerate. Finally, since $(a^{ij}_t)$ is real, its eigenvalues are symmetrically distributed with multiplicity about the real axis and thus, for all $t>0$ sufficiently small:
$$
\opSgn(a^{ij}_t) = \opSgn(a^{ij}_4) = (-1)^{n+1}\opSgn(\opHess(R_f)(p)).
$$
\noindent Finally, by \cite{Kato} (c.f. \cite{SmiDTI} for details), for $t$ small, $J_t$ has exactly one other strictly negative real eigenvalue, which is the perturbation of the eigenvalue $(-1)$ of $J_0$. Thus:
$$
\opSgn(J_t) = -(-1)^{n+1}\opSgn(\opHess(R_f)(p)).
$$
\noindent This completes the proof.\qed
\newhead{Uniqueness}
\newsubhead{Overview}
\noindent We prove under fairly weak hypotheses that the only locally strictly convex solutions are Ye's spheres.
\proclaim{Theorem \procref{ThmUniqueness}}
\noindent Suppose $K=H$ is mean curvature, or $K$ is convex and satisfies the Unbounded Growth Axiom (Axiom $(vii)$), then, for all $f\in C^\infty(M)$ and for all $C>0$, there exists $B>0$ such that if $H\geqslant B$, and if $\Sigma^n$ is a locally strictly convex immersed sphere in $M$ such that:
\medskip
\myitem{(i)} $\Sigma^n$ is of prescribed $K$-curvature equal to $H(1+H^{-2}f)$; and
\medskip
\myitem{(ii)} $\opDiam_M(\Sigma^n)\leqslant CH^{-1}$.
\medskip
\noindent Then $\Sigma^n$ is one Ye's spheres.
\endproclaim
\newsubhead{General Uniqueness}
\noindent We prove Theorem \procref{ThmUniqueness} using the following slightly weaker result: let $(k_n)_\ninn\in]0,\infty[$ be a sequence of positive real numbers converging to $+\infty$, and let $(\Sigma_n)_\ninn$ be a sequence of locally strictly convex immersed spheres in $M$ such that, for all $n$, $\Sigma_n$ has prescribed $K$-curvature equal to $k_n(1+k_n^{-2}f)$. Denote by $g$ the metric over $M$, and, for all $n$, consider the pointed Riemannian manifold $M_n:=(M,k_n^{-1}g,p_n)$. Trivally, $(M_n)_\ninn$ converges to $(\Bbb{R}^{n+1},0,g_\opEuc)$ in the $C^\infty$-Cheeger/Gromov sense. For all $n$, we now view $\Sigma_n$ as a locally strictly convex embedded submanifold of $M_n$ of constant $K$-curvature equal to $1$. We obtain:
\proclaim{Theorem \nextprocno}
\noindent If $K$ is convex, or $K=H$ is mean curvature, and if $(\Sigma_n)_\ninn$ converges in the $C^\infty$-sense to a compact, locally strictly convex immersed submanifold of $\Bbb{R}^{n+1}$, then, for sufficiently large $n$, $\Sigma_n$ is one of Ye's spheres.
\endproclaim
\proclabel{ThmWeakerUniqueness}
\noindent Let us denote the limit by $\Sigma_0$. Trivially, $\Sigma_0$ is locally convex and of constant $K$-curvature equal to $1$. It follows by Theorem \procref{ThmAlexander}, that $\Sigma_0$ is a sphere of unit radius, and, in particular, for sufficiently large $n$, $\Sigma_n$ is embedded and bounds a convex set, $\Omega_n$, say. For all $n$, we now choose $p_n\in\Omega_n$ to be any point maximising distance to $\Sigma_n$. Then $\Sigma_0$ is the unit sphere, centred on the origin, and so for sufficiently large $n$ there exists a unique function $\varphi_n\in C^\infty(S)$ such that $\Sigma_n$ is the graph of $\varphi_n$ over $S^{t_n,0}$ where $t_n:=k_n^{-1}$. We first arrange for $\varphi_n$ to lie in $H_1^\perp$.
\proclaim{Proposition \nextprocno}
\noindent For sufficiently large $n$, there exists a unique pair $(x_n,\psi_n)\in M\times H_1^\perp$ such that $\Sigma_n$ is the graph of $\psi_n$ over $S^{t_n,x_n}$.
\endproclaim
\proof We consider the mapping $\Phi:[0,\infty[\times T_pM\times C^\infty(S)\rightarrow C^\infty(S)$ defined such that:
$$
\opExp_p(t\Phi(t,V,\psi)(x)x) = \opExp_{tV}(t\psi(x)x_{tV}),
$$
\noindent where $x_{tV}$ is the parallel transport of $x$ along the unique geodesic joining $p$ to $tV$ and we identify $tV$ with the point $\opExp(tV)\in M$. $\Phi$ defines a smooth map between Frechet Spaces. Let $D_t\Phi$, $D_V\Phi$ and $D_\psi\Phi$ denote the partial derivatives of $\Phi$ with respect to the three respective variables. At $(0,0,0)$, we readily obtain:
$$\matrix
D_\psi\Phi\cdot f \hfill&= f,\hfill\cr
D_V\Phi\cdot U \hfill&= \langle U,\msf{N}\rangle,\hfill\cr
\endmatrix$$
\noindent where $\msf{N}$ is the outward pointing unit normal over $S$. It follows that $D_V\Phi\oplus D_\psi\Phi:T_pM\times H_1^\perp\rightarrow C^\infty(S)$ defines an isomorphism of Frechet spaces, and thus the restriction of $\Phi$ to $\left\{t\right\}\times T_pM\times H_1^\perp$ is invertible for all sufficiently small $t$. Thus, for $n$ sufficiently large, there exists $(V_n,\psi_n)\in T_pM\times H_1^\perp$ such that:
$$
\Phi(t_n,V_n,\psi_n) = \varphi_n.
$$
\noindent This completes the proof.\qed
\medskip
\noindent Thus, without loss of generality, we may assume that $\varphi_n\in H_1^\perp$. For all $n$, denote $\psi_n=t_n^{-2}\varphi_n$.
\proclaim{Proposition \nextprocno}
\noindent There exists $\varphi_0\in C^\infty(S)$ towards which $(\psi_n)_\ninn$ subconverges in the $C^\infty$ sense.
\endproclaim
\proclabel{PropFirstRefinement}
\proof For all $n$, the $K$-curvature of $\Sigma_n$ is prescribed by $t_n^{-1}(1+t_n^2f)$. Thus, by Proposition \procref{PropExpansionOfPertCurvature}, for all $n$:
$$
\frac{1}{n}(n+\Delta^0)\varphi_n = O(t_n^2) + O(\varphi_n^2).
$$
\noindent For all $(k,\alpha)$, we consider the $C^{k,\alpha}$-norm of $\varphi_n$. Since $(n+\Delta^0)$ defines a Frechet Space isomorphism from $H_1^\perp\minter C^\infty(S)$ to itself, there exists $K>0$ such that, for all $n$:
$$
\|\varphi_n\|_{k,\alpha} \leqslant O(t_n^2) + O(\|\varphi_n\|_{k,\alpha}^2).
$$
\noindent However, since $(\varphi_n)_\ninn$ converges to $0$ in the $C^\infty$ sense, in particular:
$$
(\|\varphi_n\|_{k,\alpha})_\ninn\rightarrow 0.
$$
\noindent There thus exists $K>0$ such that, for sufficiently large $n$:
$$
\|\varphi_n\|_{k,\alpha} \leqslant Kt_n^2.
$$
\noindent The result now follows by the Arzela-Ascoli Theorem.\qed
\medskip
\noindent We define $\varphi_{0,x_n}$ as in Proposition \procref{PropFirstAndSecondOrderTerms}, and, for all $n$, we now define $\psi_n$ by:
$$
\psi_n = t_n^{-1}(t_n^{-2}\varphi_n - \varphi_{0,x_n}).
$$
\noindent We obtain:
\proclaim{Proposition \nextprocno}
\myitem{(i)} $(x_n)_\ninn$ subconverges to a critical point of $R + \frac{2(n+3)}{(n+1)}f$; and
\medskip
\myitem{(ii)} there exists $\varphi_0\in C^\infty(S)$ towards which $(\psi_n)_\ninn$ subconverges in the $C^\infty$ sense.
\endproclaim
\proclabel{PropSecondRefinement}
\proof By compactness, we may assume that there exists $x_0\in M$ towards which $(x_n)_\ninn$ subconverges. By Proposition \procref{PropExpansionOfPertCurvature}, for all $n$:
$$
\frac{1}{n}(n+\Delta^0)\psi_n = -\frac{1}{4}\opRic^{x_n}_{yy;y} - f^{x_n}_{;y} + O(t_n^4).
$$
\noindent Composing with $\Pi_1$ yields:
$$
\frac{(n+1)}{2(n+3)}R^{x_n}_{;y} + f^{x_n}_{;y} = O(t_n^4).
$$
\noindent Taking limits yields:
$$
\frac{(n+1)}{2(n+3)}R^{x_0}_{;y} + f^{x_0}_{;y} = 0.
$$
\noindent $x_0$ is therefore a critical point of $R + \frac{2(n+3)}{(n+1)}f$, and the first assertion follows. It now follows as in the proof of Proposition \procref{PropFirstRefinement} that $(\psi_n)_\ninn$ is uniformly bounded in the $C^{k,\alpha}$-norm for all $(k,\alpha)$, and the second assertion now follows as before.\qed
\medskip
\noindent Finally, we define $\varphi_{1,x_n}$ as in Proposition \procref{PropFirstAndSecondOrderTerms}, and, for all $n$, we now define $\psi_n$ by:
$$
\psi_n = t_n^{-2}(t_n^{-2}\varphi_n - \varphi_{0,x_n} - t\varphi_{1,x_n})
$$
\noindent In like manner to the two preceeding propositions, we obtain:
\proclaim{Proposition \nextprocno}
\noindent There exists $\psi_0\in C^\infty(S)$ towards which $(\psi_n)_\ninn$ subconverges in the $C^\infty$ sense.
\endproclaim
\proclabel{PropThirdRefinement}
\noindent We may now prove Theorem \procref{ThmWeakerUniqueness}:
\medskip
{\bf\noindent Proof of Theorem \procref{ThmWeakerUniqueness}:\ }This follows immediately from Corollary \procref{CorUniqueness} and Propositions \procref{PropSecondRefinement} and \procref{PropThirdRefinement}.\qed
\newsubhead{Convex Curvature and Mean Curvature}
\noindent It now remains to show that the hypotheses of Theorem \procref{ThmWeakerUniqueness} are satisfied when $K$ is convex or when $K=H$ is mean curvature.
\proclaim{Proposition \nextprocno}
\noindent With the notation of Theorem \procref{ThmUniqueness}, $(\Sigma_n)_\ninn$ converges in the $C^\infty$ sense to the unit sphere in $\Bbb{R}^{n+1}$ about the origin.
\endproclaim
\proclabel{PropFromWeakToStrongUniqueness}
\proof Suppose first that $K$ is convex and satisfies the Unbounded Growth Axiom. By Theorem \procref{ThmAlexander}, for all sufficiently large $n$, $\Sigma_n$ is embedded and bounds a convex set, $\Omega_n$, such that $\opDiam(\Omega_n)\leqslant C$. Thus, by compactness of the family of convex sets, there exists a bounded convex set $\Omega_0\subseteq\Bbb{R}^{n+1}$ towards which $(\Omega_n)_\ninn$ subconverges. We first show that $\Omega_0$ has non-trivial interior. Suppose the contrary. Then, without loss of generality, $\Omega_0\subseteq\Bbb{R}^n$. Let $S\subseteq\Bbb{R}^{n+1}$ be a sphere, with centre in $\Bbb{R}^n$ which contains $\Omega_0$ and which is tangent to $\Omega_0$ at some point. Since $K$ satisfies the Unbounded Growth Axiom (Axiom $(vii)$), contracting $S$ in the $(n+1)$'th direction yields a hypersurface $S'$ of $K$-curvature as large as we wish at the point of contact. This is absurd, since $\Omega_0$ has constant $K$-curvature equal to $1$ in the viscosity sense (being a limit of hypersurfaces of prescribed $K$-curvature).
\medskip
\noindent Now, choose $p\in\partial\Omega_0$. If $\Omega_0$ does not satisfy the Local Geodesic Property at $p$, then by \cite{ShengUrbasWang}, $\partial\Omega_0$ is smooth in a neighbourhood of $p$. Moreover, the shape operator of $\partial\Omega_0$ is strictly positive definite at $p$, and so $\partial\Omega_0$ is also strictly convex in a neighbourhood of $p$, and, in particular, does not satisfy the Local Geodesic Property over this neighbourhood. The subset $X\subseteq\partial\Omega_0$ of points where $\partial\Omega_0$ does not satisfy the Local Geodesic Property is therefore open. However, $\Omega_0$ is bounded, and it follows from Lemma \procref{LemmaLocalGeodesicAlternative} that $X=\emptyset$. We deduce once again by \cite{ShengUrbasWang} that $\partial\Omega_0$ is smooth and that $(\Sigma_n)_\ninn$ converges in the $C^\infty$ sense to $\Sigma_0:=\partial\Omega_0$.
\medskip
\noindent If $K=H$ is mean curvature, then, by convexity, the shape operator of every $\Sigma_n$ is uniformly bounded, and it follows from the Arzela-Ascoli Theorem for immersed submanifolds (c.f. \cite{SmiAAT}) and elliptic regularity that there exists a locally convex $C^\infty$ immersed hypersurface $\Sigma_0$ towards which $(\Sigma_n)_\ninn$ subconverges in the $C^\infty$ sense modulo reparametrisation. In both cases, $\Sigma_0$ has constant $K$-curvature equal to $1$, and it therefore follows by Theorem \procref{ThmAlexandrov} that $\Sigma_0$ is a sphere of radius $1$. Since $0$ is the limit of points in $\Omega_n$ maximising distance to $\Sigma_n$, it is also the centre of $\Sigma_0$, and this completes the proof.\qed
\medskip
\noindent Theorem \procref{ThmUniqueness} follows immediately:
\medskip
{\noindent\bf Proof of Theorem \procref{ThmUniqueness}:\ }This is the union of Theorem \procref{ThmWeakerUniqueness} and Proposition \procref{PropFromWeakToStrongUniqueness}.\qed
\newhead{Bibliography}
{\leftskip = 5ex \parindent = -5ex
\leavevmode\hbox to 4ex{\hfil \cite{Alexander}}\hskip 1ex{Alexander S., Locally convex hypersurfaces of negatively curved spaces, {\sl Proc. Amer. Math. Soc.} {\bf 64} (1977), no. 2, 321--325}
\medskip
\leavevmode\hbox to 4ex{\hfil \cite{Almgren}}\hskip 1ex{Almgren F. J. Jr., {\sl Plateau's problem: An invitation to varifold geometry},.W. A. Benjamin, Inc., New York-Amsterdam, (1966)}
\medskip
%\leavevmode\hbox to 4ex{\hfil \cite{Berger}}\hskip 1ex{Berger M., {\sl A panoramic view of Riemannian geometry}, Springer-Verlag, Berlin, (2003)}
%\medskip
%\leavevmode\hbox to 4ex{\hfil \cite{BoehmWilking}}\hskip 1ex{B\"ohm C., Wilking B., Manifolds with positive curvature operators are space forms, {\sl Ann. of Math. (2)} {\bf 167} (2008), no. 3, 1079--1097}
%\medskip
\leavevmode\hbox to 4ex{\hfil \cite{CaffNirSprV}}\hskip 1ex{Caffarelli L., Nirenberg L., Spruck J., Nonlinear second-order elliptic equations. V. The Dirichlet problem for Weingarten hypersurfaces, {\sl Comm. Pure Appl. Math.} {\bf 41} (1988), no. 1, 47--70}
\medskip
\leavevmode\hbox to 4ex{\hfil \cite{Dold}}\hskip 1ex{Dold A., {\sl Lectures on Algebraic Topology}, Classics in Mathematics, Springer-Verlag, Berlin, Heidelberg, (1995)}
\medskip
\leavevmode\hbox to 4ex{\hfil \cite{ElwTrombI}}\hskip 1ex{Elworthy K. D., Tromba A. J., Differential structures and Fredholm maps on Banach manifolds, in {\sl Global Analysis (Proc. Sympos. Pure Math.)}, Vol. {\bf XV}, Berkeley, Calif., (1968), 45--94 Amer. Math. Soc., Providence, R.I.}
\medskip
\leavevmode\hbox to 4ex{\hfil \cite{ElwTrombII}}\hskip 1ex{Elworthy K. D., Tromba A. J., Degree theory on Banach manifolds, in {\sl Nonlinear Functional Analysis (Proc. Sympos. Pure Math.)}, Vol. {\bf XVIII}, Part 1, Chicago, Ill., (1968), 86--94, Amer. Math. Soc., Providence, R.I.}
\medskip
\leavevmode\hbox to 4ex{\hfil \cite{EspGalMira}}\hskip 1ex{Espinar J. M., G\'alvez J. A., Mira P., Hypersurfaces in $\Bbb H^{n+1}$ and conformally invariant equations: the generalized Christoffel and Nirenberg problems, {\sl J. Eur. Math. Soc.}, {\bf 11}, no. 4, (2009), 903--939}
\medskip
\leavevmode\hbox to 4ex{\hfil \cite{RosEsp}}\hskip 1ex{Espinar J., Rosenberg H., When strictly locally convex hypersurfaces are embedded, arXiv:1003.0101}
\medskip
\leavevmode\hbox to 4ex{\hfil \cite{Fethi}}\hskip 1ex{Fethi M., Constant $k$-curvature hypersurfaces in Riemannian manifolds, {\sl Differential Geom. Appl.} {\bf 28} (2010), no. 1, 1--11}
\medskip
\leavevmode\hbox to 4ex{\hfil \cite{Kato}}\hskip 1ex{Kato T., {\sl Perturbation theory for linear operators}, Die Grundlehren der mathematischen Wissenschaften, {\bf 132}, Springer-Verlag, New York, (1966)}
\medskip
\leavevmode\hbox to 4ex{\hfil \cite{LabB}}\hskip 1ex{Labourie F., Probl\`eme de Minkowski et surfaces \`a courbure constante dans les vari\'et\'es hyperboliques, {\sl Bull. Soc. Math. France}, {\bf 119}, no. 3,(1991), 307--325}
\medskip
\leavevmode\hbox to 4ex{\hfil \cite{LabI}}\hskip 1ex{Labourie F., Immersions isométriques elliptiques et courbes pseudoholomorphes, {\sl J. Diff. Geom.} {\bf 30}, (1989), 395-44}
\medskip
\leavevmode\hbox to 4ex{\hfil \cite{PacardXu}}\hskip 1ex{Pacard F., Xu X., Constant mean curvature spheres in Riemannian manifolds,\break {\sl Manuscripta Math.} {\bf 128} (2009), no. 3, 275--295}
\medskip
%\leavevmode\hbox to 4ex{\hfil \cite{MicallefMoore}}\hskip 1ex{Micallef M. J., Moore J. D., Minimal two-spheres and the topology of manifolds with positive curvature on totally isotropic two-planes, {\sl Ann. of Math.} (2), \bf{127}, (1988), no. 1, 199--227}
%\medskip
\leavevmode\hbox to 4ex{\hfil \cite{Nirenberg}}\hskip 1ex{Nirenberg L., The Weyl and Minkowski problem in differential geometry in the large, {\sl Comm. Pure Appl. Math.}, {\bf 6}, no. 3, (1953), pp. 337--394}
\medskip
\leavevmode\hbox to 4ex{\hfil \cite{Pogorelov}}\hskip 1ex{Pogorelov A.V., Extrinsic geometry of convex surfaces, {\sl Israel program for scientific translation}, Jerusalem, (1973)}
\medskip
\leavevmode\hbox to 4ex{\hfil \cite{SmiDTI}}\hskip 1ex{Rosenberg H., Smith G., Degree Theory of Immersed Hypersurfaces, in preparation}
\medskip
\leavevmode\hbox to 4ex{\hfil \cite{SchneiderI}}\hskip 1ex{Schneider M., Alexandrov Embedded closed magnetic geodesics on $S^2$,\break arXiv:0903.1128}
\medskip
\leavevmode\hbox to 4ex{\hfil \cite{SchneiderII}}\hskip 1ex{Schneider M., Closed magnetic geodesics on closed hyperbolic Riemann Surfaces, arXiv:1009.1723}
\medskip
%\leavevmode\hbox to 4ex{\hfil \cite{Schwarz}}\hskip 1ex{Schwarz M., {\sl Morse Homology}, Progress in Mathematics, %{\bf 111}, Birkh\"auser Verlag, Basel, Boston, Berlin, (1993)}
%\medskip
\leavevmode\hbox to 4ex{\hfil \cite{ShengUrbasWang}}\hskip 1ex{Sheng W., Urbas J., Wang X., Interior curvature bounds for a class of curvature equations. (English summary), {\sl Duke Math. J.} {\bf 123} (2004), no. 2, 235--264}
\medskip
\leavevmode\hbox to 4ex{\hfil \cite{SmiAAT}}\hskip 1ex{Smith G., An Arzela-Ascoli Theorem for Immersed Submanifolds, {\sl Ann. Fac. Sci. Toulouse Math.} {\bf 16} (2007), no. 4, 817--866}
\medskip
\leavevmode\hbox to 4ex{\hfil \cite{SmiSLC}}\hskip 1ex{Smith G., Special Lagrangian curvature, to appear in {\sl Math. Annalen}}
\medskip
%\leavevmode\hbox to 4ex{\hfil \cite{SmiCGC}}\hskip 1ex{Smith G., Constant Gaussian Curvature hypersurfaces in Hadamard manifolds,\break arXiv:0908.3590}
%\medskip
\leavevmode\hbox to 4ex{\hfil \cite{SmiNLD}}\hskip 1ex{Smith G., The non-linear Dirichlet problem in Hadamard manifolds, arXiv:0908.3590}
\medskip
\leavevmode\hbox to 4ex{\hfil \cite{SmiPPG}}\hskip 1ex{Smith G., The Plateau problem for general curvature functions, arXiv:1008.3545}
\medskip
\leavevmode\hbox to 4ex{\hfil \cite{Stein}}\hskip 1ex{Stein E., {\sl Singular integrals and differentiability properties of functions}, Princeton University Press, (1970)}
\medskip
\leavevmode\hbox to 4ex{\hfil \cite{Tromb}}\hskip 1ex{Tromba A. J., The Euler characteristic of vector fields on Banach manifolds and a globalization of Leray-Schauder degree, {\sl Adv. in Math.}, {\bf 28}, (1978), no. 2, 148--173}
\medskip
\leavevmode\hbox to 4ex{\hfil \cite{White}}\hskip 1ex{White B., The space of $m$-dimensional surfaces that are stationary for a parametric elliptic functional, {\sl Indiana Univ. Math. J.} {\bf 36}, (no. 3), (1987), 567--602}
\medskip
\leavevmode\hbox to 4ex{\hfil \cite{Ye}}\hskip 1ex{Ye R., Foliation by constant mean curvature spheres, {\sl Pacific J. Math.} {bf 147} (1991), no. 2, 381--396}
\par
}%
\enddocument